%% file: manuscript_revision__1_.tex
\documentclass[onefignum,onetabnum]{siamart171218}

\input{ex_shared}
\usepackage{subfig}
\ifpdf
\hypersetup{
  pdftitle={Boundary-integral simulations of open membranes},
  pdfauthor={Han Zhou, Y.-N. Young, and Yoichioro Mori}
}
\fi

\myexternaldocument{ex_supplement}

\DeclareMathOperator{\Id}{Id}

\newcommand{\eqdef}{\overset{\mbox{\tiny{def}}}{=}}

\newcommand{\paren}[1]{\left(#1\right)}

\newcommand{\jump}[1]{[\![#1]\!]}

\newcommand{\at}[2]{\left. #1 \right|_{#2}}
\newcommand{\grad}[1]{\nabla #1}

\newcommand{\mc}[1]{\mathcal{#1}}

\newcommand{\wt}[1]{\widetilde{#1}}
\newcommand{\bm}[1]{\boldsymbol{#1}}

\newcommand{\norm}[1]{\left\lVert #1 \right\rVert}
\newcommand{\dual}[2]{\left\langle #1,#2 \right\rangle}

\newcommand{\vph}{\varphi}

\def\tsc#1{\csdef{#1}{\textsc{\lowercase{#1}}\xspace}}
\tsc{WGM}
\tsc{QE}
\tsc{EP}
\tsc{PMS}
\tsc{BEC}
\tsc{DE}

\begin{document}

\maketitle

\begin{abstract}
  In this work, we present a mathematical and computational framework to model the dynamics of open lipid bilayer membranes interacting with ambient Stokes flow. The model explicitly couples the three-dimensional viscous fluid, the two-dimensional membrane surface, and its one-dimensional free edge. We develop an axisymmetric hybrid BEM-FEM method that solves the problem with an effective one-dimensional formulation. A key component is a local mesh refinement strategy designed to accurately resolve singularities and boundary layers originating at the membrane edge. Several numerical examples are provided to showcase its ability to capture intricate edge dynamics and multiscale fluid-membrane coupling.
\end{abstract}



\section{Introduction}

Lipid bilayer membranes are fundamental structural components of cellular organelles, as they self-enclose to form a vesicle in a viscous fluid to separate interior from exterior contents. As a biomimetic model for blood cells, the hydrodynamics of giant vesicles provides fundamental insights into the mechanics of biomembranes, namely the coupling between membrane and fluid mechanics \cite{Saitoh1998,Misbah2006, Cooper2015,Arroyo2009, Aubin2016,Lai2001, Cortez2005}. 
Under normal conditions, lipid bilayer membranes close on themselves in a viscous polar solvent (such as water) to form a vesicle \cite{seifert1997configurations}.
However, in various biological and biophysical processes, the lipid bilayer membrane may open or reseal subjected to a mechanical stress or a chemical gradient \cite{ryham2011aqueous,Cooper2015, Tang2017, Klenow2024}.
Various techniques are developed to control the opening and closing of membranes by applying osmotic stress, using light-sensitive molecules, or employing high pulse of electric field. 
These techniques enable precise manipulation of membrane opening and closing, with applications spanning from targeted drug delivery to synthetic vesicle engineering \cite{Malik2022, Saitoh1998, Srividya2008, Umeda2005}.

Despite their relevance and importance, open membranes pose significant numerical challenges due to the presence of free edges. Unlike closed membranes, which are continuous and without boundaries, open membranes have edges where the lipid bilayer tilts to minimize exposure to solvent molecules \cite{ryham2013teardrop,fu2020simulation}. These membrane edges require boundary conditions that enforce local balances of forces and torques along the edge \cite{Asgari2015, Srividya2008, Alexander2006}. Such conditions are typically expressed as constraints involving line tension, bending moments, and the force balance between the membrane and the surrounding fluid. Most existing studies on open membranes focus on equilibrium configurations without hydrodynamic coupling \cite{Biria2013, Laadhari2010, Capovilla2002, Deserno2015, Palmer2021, Palmer2022, Tu2010, Tu2011}. In particular, axisymmetric open membranes have been analyzed to understand their stable equilibrium shapes, the effects of edge tension, and morphological transitions \cite{Palmer2024, Tu2003, Zhou2018, Yin2005, Yin2005b}.
For open membranes with hydrodynamics, several studies have been conducted under the assumption that the membrane maintains a prescribed shape, such as a spherical cap \cite{Aubin2016, Ryham2018, Arroyo2009}. Diffuse interface methods have also been developed to model open membranes coupled with hydrodynamics \cite{Cohen2012, Cohen2014}, or under purely geometric gradient flows \cite{Wang2007}. However, the sharp-interface limits of such models remain unclear, raising questions about their convergence. Fully resolved sharp-interface formulations for the coupled dynamics of open membranes and Stokes flow remain largely unexplored.

For closed fluid–membrane interaction problems, a variety of numerical algorithms have been developed, including the immersed boundary method \cite{Kim2010, Lai2019, Chuan2021, Lai2022}, the boundary integral method \cite{Veerapaneni2009, Farutin2014, HSU2019747}, the parametric finite element method \cite{Barrett2020}, and the diffuse interface method \cite{Bottacchiari2024, Reuther2016}.
{\color{black} In addition, hybrid methods have been proposed that incorporate membrane and bulk viscosities and simulate vesicle electrohydrodynamics \cite{Boedec2017,hu2015hybrid,hu2016vesicle}. }

{\color{black}For open membranes, parametric finite-element formulations have been widely employed to simulate the gradient flow of the Helfrich energy and obtain equilibrium membrane shapes without explicit hydrodynamic coupling to the surrounding viscous fluid \cite{balchunas2020force,balchunas2019equation,Barrett2020,Barrett2017,Barrett2021}. Incorporating fluid–membrane coupling would require solving the Stokes equations in domains bounded by open and deforming surfaces, a task complicated by the singular behavior near the free edge—the so-called edge condition \cite{Hayashi1977,Kirvalidze1997,Stephan1987}. Related numerical strategies have been developed in boundary-integral formulations for Laplace and Helmholtz problems on open curves, providing conceptual and computational approaches to regularizing such edge singularities \cite{Jiang2004,Jiang2003,Helsing2024}.}


This work develops a mathematical and computational framework for simulating the dynamics of an open axisymmetric lipid bilayer membrane immersed in an incompressible Stokesian fluid. 
The system is inherently multiscale, involving the coupling of PDEs across domains of differing dimensions: three-dimensional (3D) bulk fluid regions, two-dimensional (2D) membrane surfaces, and one-dimensional (1D) membrane edges. Such mixed-dimensional PDE systems are essential for accurately capturing physical phenomena governed by cross-dimensional interactions. While prior studies have addressed 3D–1D coupling in applications such as blood flow and oxygen transport in microcirculation \cite{quarteroni2007stability, kuchta2020mixed, grappein2024xfem}, models that simultaneously and consistently integrate 3D fluid dynamics, 2D membrane elasticity, and 1D edge line tension remain scarce. 
To solve the resulting mixed-dimensional system, we employ a boundary integral formulation with axisymmetric reduction and develop a hybrid boundary element–finite element method (BEM-FEM) for the effectively one-dimensional problem.
The method is further enhanced by a local mesh refinement strategy to resolve the singular behavior near the membrane edge. The proposed model and numerical scheme are validated through a series of simulations, with results benchmarked against existing analytical and computational studies. To our knowledge, this is the first fully coupled and multiscale computational framework for modeling open lipid membrane dynamics in Stokes flow.

The remainder of this paper is organized as follows.
Section~\ref{sec:2} introduces the membrane energy functional and the fluid–membrane interaction, which together yield the governing equations.
Section~\ref{sec:3} presents the axisymmetric reduction and weak formulation, reducing the problem to a one-dimensional representation.
Section~\ref{sec:4} describes the numerical method used to solve the coupled system, with particular attention to the treatment of edge singularities.
Section~\ref{sec:5} presents numerical experiments that demonstrate the effectiveness of the proposed method.
Finally, Section~\ref{sec:6} concludes with a summary of the work and outlines directions for future research.

\section{Model}\label{sec:2}
This study focuses on modeling and simulating an open inextensible elastic membrane interacting with a viscous fluid in the Stokes regime (no inertial effect). In this section, we employ an energetic variational approach to derive the governing equations for the dynamics of the membrane.

\subsection{Differential geometry of the open membrane}
\begin{figure}[htbp]
    \centering
    \includegraphics[width=0.3\linewidth]{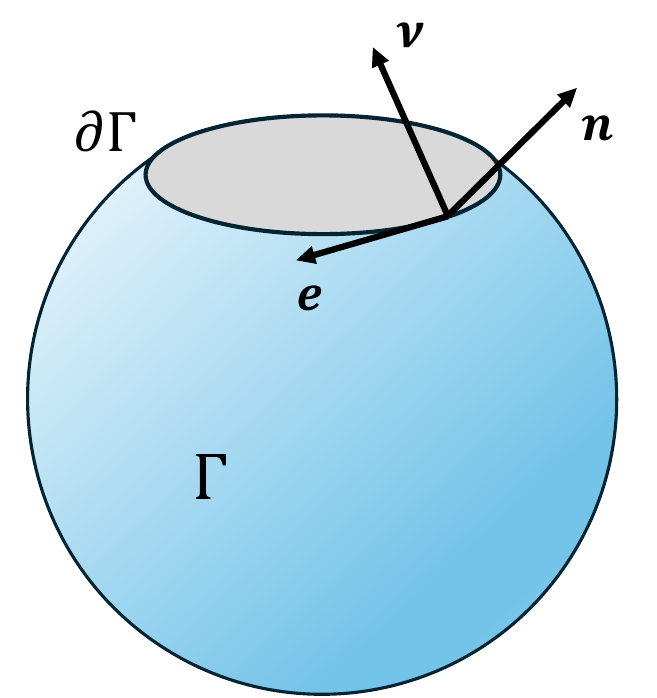}
    \caption{Illustration for an open membrane surface $\Gamma$ with boundary $\partial\Gamma$. $\bm{n}$ is the unit normal of $\Gamma$, $\bm{e}$ is the tangent and $\bm{\nu}$ is the co-normal along the edge $\partial\Gamma$.}
    \label{fig:illustration}
\end{figure}

 To present the governing equations, we first introduce some key differential geometric concepts for an open surface and its boundary. Consider an open inextensible, elastic membrane with one or multiple holes. The membrane is assumed to be a smooth and orientable open surface $\Gamma\subset \mathbb{R}^3$ with zero thickness (Figure~\ref{fig:illustration}).  The boundary $\partial\Gamma$ corresponds to the open edge of elastic membrane $\Gamma$. Suppose the membrane is immersed in an incompressible viscous fluid in the domain $\Omega = \mathbb{R}^3\setminus \overline{\Gamma}$. 

Consider the surface characterized by the parametrization $\bm X (\bm{\theta})\in \Gamma\subset\mathbb{R}^3$, where $\bm{\theta} = (\theta_1, \theta_2)$ denotes a coordinate system defined on the open set $\mathcal{U}$. The metric tensor is defined as
{\color{black}
\begin{equation}
g_{ij} = \partial_i \bm X \cdot \partial_j \bm X, \quad \text{for} \ i,j = 1,2.
\end{equation}
}
The dual tangent vectors are given by
{\color{black}
\begin{equation}
    \bm X ^{*}_{1} = g^{11} \partial_1 \bm X  + g^{12} \partial_2 \bm X  ,\quad\bm X ^{*}_{2} = g^{21} \partial_1 \bm X  + g^{22} \partial_2 \bm X  ,
\end{equation}
}
where $(g^{ij})=(g_{ij})^{-1}$ represents the inverse of the metric tensor.
Utilizing the dual tangent vectors, the surface gradient and surface divergence of a scalar function $f$ and a vector field $\bm g$ defined on $\Gamma$ can be represented as
\begin{equation}
    \grad_{\Gamma} f = \frac{\partial f}{\partial \theta_1}\bm X^*_1 +\frac{\partial f}{\partial \theta_2}\bm X^*_2 ,\quad  \grad_{\Gamma}\cdot \bm g= \frac{\partial \bm g}{\partial \theta_1} \cdot\bm X^*_1 +\frac{\partial \bm g}{\partial \theta_2}\cdot\bm X^*_2 
\end{equation}
Let {\color{black}$\bm n = |\partial_1 \bm X \times \partial_2\bm X |^{-1}(\partial_1 \bm X \times \partial_2\bm X )$} be the unit normal vector field.
The orthogonal projection onto the tangent space of $\Gamma$, denoted as $\bm{P}_{\Gamma}$, is given by
\begin{align}
\label{eq:P_gamma}
   \textstyle \bm{P}_{\Gamma} = \bm{I} - \bm{n}\otimes\bm{n}  = \sum_{i=1}^2\bm{\tau}_i\otimes\bm{\tau}_i,
\end{align}
where $\bm{\tau}_i$ denotes orthogonal unit tangent vectors on the surface.
Defining the curvature tensor $\bm{L} = -\nabla_{\Gamma}\bm{n}$ (also referred to as the Weingarten map), the mean and Gaussian curvatures $H$ and $K$ can be expressed as
\begin{equation}
    H = \frac{1}{2}\text{tr} \bm{L} , \quad K =\frac{1}{2}\left(\left(\text{tr} \bm{L}\right)^2 - \text{tr}\left(\bm{L}^2\right)\right).
\end{equation}

In Figure~\ref{fig:illustration}, the unit tangent $\bm{e}$ along the edge $\partial\Gamma$ is chosen such that {\color{black}the outward-pointing unit co-normal vector} $\bm{\nu}$ of $\partial\Gamma$, given by $\bm{\nu} = \bm{e}\times\bm{n}$, is oriented away from $\Gamma$. The Darboux frame for $\partial\Gamma$ is denoted as $\{\bm{e}, \bm{n}, \bm{\nu}\}$, and it satisfies the following relationships
\begin{equation}
    \bm{e}' = \kappa_n \bm{n} - \kappa_g \bm{\nu},\quad
    \bm{n}' = -\kappa_n \bm{e} -\tau_g \bm{\nu},\quad
    \bm{\nu}' = \kappa_g \bm{e} + \tau_g \bm{n},
\end{equation}
where the prime denotes the derivative with respect to the arc length parameter $s$ of $\partial\Gamma$, $\kappa_n$ and $\kappa_g$ represent the normal and geodesic curvatures of $\partial\Gamma$, and $\tau_g$ stands for the geodesic torsion of $\partial\Gamma$.

\subsection{Energy variational formulations}
A lipid bilayer membrane in a given configuration stores free energy, which drives the system toward configurations of lower energy by generating conservative forces. For an open lipid bilayer membrane, the total free energy consists of contributions from both the curvature elasticity of the membrane surface and the line tension of the edge. Mathematically, the total free energy $\mc E$ is expressed as
\begin{equation}
    \mathcal{E} = \mathcal{E}_{\Gamma} + \mathcal{E}_{\partial\Gamma}.
\end{equation}
The curvature elasticity of the membrane is modeled using the Helfrich energy,
\begin{equation}
    \mathcal{E}_{\Gamma} = \int_{\Gamma} \paren{\alpha \paren{H - c_0}^2 + \alpha_G K} \, dA,
\end{equation}
where $ H $ is the mean curvature, $ K $ is the Gaussian curvature, $ c_0 $ is the spontaneous curvature, and $ \alpha $ and $ \alpha_G $ are the bending rigidities associated with the mean and Gaussian curvatures, respectively. The spontaneous curvature $ c_0 $ reflects the natural geometry of the membrane due to its molecular composition or external environment. If $ c_0 = 0 $ and $ \alpha_G = 0 $, the Helfrich energy reduces to the simpler Willmore energy, which depends only on the mean curvature.  For the line tension of the edge, we consider the following energy
\begin{equation}
    \mathcal{E}_{\partial\Gamma} = \gamma \int_{\partial\Gamma} 1 \, ds,
\end{equation}
where $ \gamma $ is the constant line tension, representing the energetic cost per unit length of the edge. The line tension reflects the resistance of the edge to lengthening, arising from molecular interactions at the boundary.

The conservative forces produced by the membrane can be determined by calculating the first variation of the total energy. Let's consider an arbitrary perturbation of the surface $\Gamma$ in the direction of $\bm Y $. We compute the first variation of the total energy by
\begin{equation}\label{eqn:dEdx}
\begin{aligned}
&\delta_{\bm X }\mathcal{E}[\bm{X}](\bm{Y}) = \at{\frac{d}{d\varepsilon} \mathcal{E}[\bm{X}+\varepsilon\bm{Y}]}{\varepsilon=0}\\
&=\int_{\Gamma} \alpha\left( \Delta_{\Gamma} H + 2(H-c_0)(H^2+Hc_0-K) \right) (\bm{n} \cdot \wt{\bm Y} )\,dA \\
&\, + \int_{\partial\Gamma} \left( \alpha(H-c_0)^2 + \alpha_G K + \gamma \kappa_g \right)(\bm{\nu}\cdot\wt{\bm Y}) - (\alpha\nabla_{\Gamma}H\cdot\bm{\nu} - \alpha_G \tau_g' + \gamma \kappa_n)(\bm{n}\cdot\wt{\bm Y} ) \\
&\, + (\alpha(H-c_0) + \alpha_G \kappa_n) \bm{\nu} \cdot\nabla_{\Gamma} (\bm n\cdot\wt{\bm Y} )\,ds,
\end{aligned}
\end{equation}
where $\wt{\bm Y} = \bm Y\circ\bm X^{-1}$.
The velocities of the bulk fluid flow and the surface flow are denoted by $\bm{u}$ and $\bm{U}$, respectively. Under the assumptions of an impermeable membrane with a no-slip condition, the compatibility constraint $\bm u = \bm U$ is imposed on the surface $\Gamma$. The membrane's evolution is governed by the geometric evolution equation given by
\begin{equation}\label{eqn:kinematic}
\partial_t \bm X (\bm \theta, t) = \bm U (\bm X (\bm \theta, t), t) = \bm u (\bm X (\bm \theta, t), t), \quad (\bm \theta, t)\in\mc U\times[0,T].
\end{equation}
Furthermore, it is assumed that the fluid approaches rest at infinity, adhering to the condition
\begin{equation}
|\bm{u}(\bm x)| \to 0 \quad \text{as } |\bm{x}| \to \infty.
\end{equation}

The fluid environment within cells is commonly described as a low Reynolds number flow \cite{Arroyo2009, Malik2022}, typically modeled using the Stokes equation. Additionally, it has been proposed that dissipation resulting from membrane viscosity plays a crucial role for large vesicles and introduces a significant regularization effect for small pores \cite{Arroyo2009, Biria2013, Jia2022}. We consider both the bulk viscosity and membrane viscosity in formulating the model. Let $\mu$ and $\mu_{\Gamma}$ represent the constant viscosity coefficients of the bulk fluid and the membrane, respectively. The viscous dissipation in the bulk fluid and on the membrane surface can be characterized by the Rayleigh dissipation functions defined as
\begin{align}
    \mathcal{R}[\bm{u}] = \int_{\Omega} \mu \bm{D}(\bm{u}) \colon \bm{D}(\bm{u}) \, d\bm{x},\quad
    \mathcal{R}_{\Gamma}[\bm{U}] = \int_{\Gamma} \mu_{\Gamma} \bm{D}_{\Gamma}(\bm{U}) \colon \bm{D}_{\Gamma}(\bm{U}) \, dA,
\end{align}
where $\bm{D}(\bm{u})$ and $\bm{D}_{\Gamma}(\bm{U})$ are the rate of deformation tensors for the bulk and surface fluid flows, defined as
\begin{align}
    \bm{D}(\bm{u}) = \frac{1}{2}\paren{ \nabla\bm{u} + \nabla\bm{u}^{\top}},\quad
    \bm{D}_{\Gamma}(\bm{U}) = \frac{1}{2}\bm{P}_{\Gamma}\paren{\nabla_{\Gamma}\bm{U} + \nabla_{\Gamma}\bm{U}^{\top}}\bm{P}_{\Gamma}.
\end{align}
Here, $\bm{P}_{\Gamma}$ is the orthogonal projection operator defined in Equation~(\ref{eq:P_gamma}). 
{\color{black}The membrane deformation-rate tensor $\bm D_\Gamma(\bm U)$ is the tangential projection of the rate of change of the surface metric tensor. It depends not only on the tangential velocity but also on the normal velocity, owing to the curvature of the surface.}
The fluid dynamics under consideration pertain to incompressible bulk Stokes flow where the bulk velocity field is divergence-free, i.e., $\nabla \cdot \bm{u} = 0$.
To account for the incompressibility constraint in the bulk, we construct a Lagrangian functional related to the viscous dissipation within the bulk based on the Rayleigh dissipation function, 
\begin{equation}
    L[\bm{u}, p] = \int_{\Omega} \mu \bm{D}(\bm{u}) \colon \bm{D}(\bm{u}) - p \left(\nabla \cdot \bm{u}\right) \, d\bm{x},
\end{equation}
where the pressure $p$ in the bulk is a Lagrange multiplier for enforcing the incompressibility condition.
Moreover, the assumption of local inextensibility of the lipid bilayer membrane implies $
\nabla_{\Gamma} \cdot \bm{U} = 0$.
To address this constraint, a Lagrangian functional associated with the viscous dissipation on the membrane surface is constructed as
\begin{equation}
    L_{\Gamma}[\bm{U}, P] = \int_{\Gamma} \mu_{\Gamma} \bm{D}_{\Gamma}(\bm{U}) \colon \bm{D}_{\Gamma}(\bm{U}) - P \left(\nabla_{\Gamma} \cdot \bm{U}\right) \, dA,
\end{equation}
where the Lagrange multiplier $P$ enforces local area incompressibility and corresponds to the negative surface tension. The dissipative force can be determined by calculating the first variation of the Lagrangian functionals. Let $\bm v$ and $\bm V$ denote arbitrary perturbations of the bulk and surface velocities, respectively. As $\bm{x}$ tends to infinity, the perturbation decays such that $|\bm{v}|\to 0$ as $|\bm{x}|\to\infty$. The first variations of $L$ can be expressed as follows
\begin{equation}\label{eqn:dLdu}
\begin{aligned}
    \delta_{\bm{u}} L[\bm{u},p](\bm{v}) &= \at{\frac{d}{d\varepsilon} L[\bm{u}+\varepsilon\bm{v}, p]}{\varepsilon=0} 
    = \int_{\Omega} \paren{2\mu \bm{D}(\bm{u}) - p\bm{I}}\colon\bm{D}(\bm{v})\,d\bm{x} \\&
    = \int_{\Omega} \paren{-\nabla\cdot \bm{\sigma} }\bm{v}\,d\bm{x} + \int_{\Gamma}\jump{\bm{\sigma}\bm{n}} \cdot\bm{v} \,dA,
\end{aligned}
\end{equation}
where the stress tensor $\bm{\sigma}$ is defined as $\bm{\sigma}\paren{\bm{u}, p} = -p\bm{I} + 2 \mu \bm{D}(\bm{u})$,
and $\jump{\bm{\sigma}\bm{n}}$ denotes the jump of normal stress on $\Gamma$, given by
\begin{equation}
    \jump{\bm{\sigma}\bm{n}} (\bm{x})= \lim_{\varepsilon\to 0} \paren{\bm{\sigma}(\bm{x}-\varepsilon\bm{n}(\bm{x})) - \bm{\sigma}(\bm{x}+\varepsilon\bm{n}(\bm{x})) }\bm{n}(\bm{x}), \quad\bm{x}\in\Gamma.
\end{equation}
Similarly, the first variations of $L_{\Gamma}$ can be expressed as
\begin{equation}\label{eqn:dLgdU}
\begin{aligned}
    \delta_{\bm{U}}L_{\Gamma}[\bm{U},P](\bm{V}) &= \at{\frac{d}{d\varepsilon} L_{\Gamma}[\bm{U}+\varepsilon\bm{V}, P]}{\varepsilon=0} 
    = \int_{\Gamma} \paren{2\mu_{\Gamma}\bm{D}_{\Gamma}(\bm{U}) - P\bm{I}} \colon \bm{D}_{\Gamma}(\bm{V}) \,dA 
    \\& = \int_{\Gamma} \paren{-\nabla_{\Gamma}\cdot \bm{\sigma}_{\Gamma}} \cdot \bm{V} \, dA + \int_{\partial\Gamma} \paren{\bm{\sigma}_{\Gamma}\bm{\nu}}\cdot\bm{V}\,ds,
\end{aligned}
\end{equation}
where {\color{black}the stress tensor for the membrane fluid flow $\bm{\sigma}_{\Gamma}$} is defined as $\bm{\sigma}_{\Gamma}\paren{\bm{U}, P} = -P\bm{P}_{\Gamma} +2 \mu_{\Gamma} \bm{D}_{\Gamma}(\bm{U})$.

\subsection{Governing equations}
The governing equations are derived via the principle of maximum dissipation, formulated as a variational condition on the rate of energy dissipation (see for example \cite{hyon2010energetic}). For our problem, this can be written as follows
\begin{equation}
\delta_{\bm{X}}\mathcal{E}[\bm{X}](\bm{Y})=-\paren{\delta L[\bm{u},p](\bm{v})+\delta L_\Gamma[\bm{U},P](\bm{V})}, 
\end{equation}
where the three terms above were computed in \eqref{eqn:dEdx}, \eqref{eqn:dLdu}, and \eqref{eqn:dLgdU} respectively. The above is to be satisfied for any arbitrary variations $\bm{Y},\bm{V}$ and $\bm{v}$ that satisfy the kinematic constraint $\bm{Y}\circ\bm X^{-1}=\bm{V}=\bm{v}$ on $\Gamma$ in Equation~\eqref{eqn:kinematic}. 

Now we summarize the governing equations for an open membrane immersed in a viscous fluid. The fluid flow in the bulk is described by the Stokes equation
\begin{equation}
    -\mu \Delta\bm u + \grad p = \bm{0},\quad \nabla\cdot\bm{u} = 0,\quad \text{ in } \Omega,
\end{equation}
with interfacial boundary condition states that the bulk fluid stress jump is balanced by bending forces and surface viscous forces on the membrane $\Gamma$, given by,
\begin{equation}
    \nabla_{\Gamma}\cdot\bm{\sigma}_{\Gamma} = \jump{\bm{\sigma}\bm{n}}+ \alpha\paren{ \Delta_{\Gamma} H + 2\paren{H-c_0} \paren{H^2+Hc_0-K}}\bm{n} ,\quad \nabla_{\Gamma} \cdot \bm{U} = \bm{0},  \quad \text{ on }\Gamma, \label{eqn:sf-force}
\end{equation}
supplemented by boundary conditions at the edge $\partial\Gamma$,
\begin{subequations}
    \begin{align}
        \bm{\sigma}_{\Gamma}\bm{\nu} +  \paren{\alpha\paren{H-c_0}^2 + \alpha_G K + \gamma \kappa_g } \bm{\nu}  &= \bm{0},\quad \text{ on }\partial\Gamma,\label{eqn:bc-conormal}\\
    \alpha\nabla_{\Gamma}H\cdot\bm{\nu} - \alpha_G \tau_g' + \gamma \kappa_n &= 0, \quad \text{ on }\partial\Gamma,\label{eqn:bc-normal}\\
    \alpha\paren{H-c_0} + \alpha_G \kappa_n &= 0,\quad \text{ on }\partial\Gamma.\label{eqn:bc-moment}
    \end{align}
\end{subequations}
Equations \eqref{eqn:bc-conormal} and \eqref{eqn:bc-normal} express the force balance conditions in the co-normal and normal directions, respectively. Equation \eqref{eqn:bc-moment} expresses the torque-free condition at $\partial\Gamma$. It is worth mentioning that our formulation captures not only the interaction of the bulk fluid with the membrane surface but also incorporates the membrane edge, {\color{black}a codimension-two object} immersed in 3D space. The forces acting on the edge are transmitted indirectly to the bulk fluid through the viscous surface,  thereby circumventing the difficulties typically associated with modeling slender-body interactions.

Our derivation ensures that in the above system the rate of change in energy is balanced by energy dissipation. Indeed, it can be shown that the above equations yield the following identity
\begin{equation}
    \frac{d}{dt}\mathcal{E} = -2 \int_{\Omega} \mu|\bm{D}(\bm{u})|^2 \,d\bm{x} -2  \int_{\Gamma}\mu_{\Gamma} |\bm{D}_{\Gamma}(\bm{U})|^2 \, dA.
\end{equation}
where $|\bm{D}(\bm{u})|$ and $|\bm{D}_{\Gamma}(\bm{U})|$ are the Frobenius norms of the deformation tensors. 

\section{Axisymmetric and weak formulation}\label{sec:3}
The governing equations for an open membrane in Stokes flow involve coupling the surface evolution equations with the bulk Stokes equation. 
We employ boundary integral formulations for the bulk Stokes equation, reducing the bulk-surface coupled PDE system to a set of differential-integral equations confined to the surface.

\subsection{Boundary integral and weak form}
Let $G(\bm{x})$ be the Stokeslet in 3D free spaces, whose entries are given by
\begin{equation}\label{eqn:stokeslet}
    G_{ij}(\bm{x}) = -\frac{1}{8\pi}\paren{\frac{\delta_{ij}}{|\bm{x}|} + \frac{x_i x_j}{|\bm{x}|^3}}.
\end{equation}
For a surface density function $\bm{f}$, define the single-layer integral operator $\mathcal{S}$ as
\begin{equation}
     (\mathcal{S}[\bm{f}])_i(\bm{x}) = - \frac{1}{\mu}\int_{\Gamma} G_{ij}\paren{\bm{x} - \bm{x}'} f_j\paren{\bm{x}'} \,d\mu(\bm{x}'), \quad \bm{x}\in\mathbb{R}^3,
\end{equation}
where we have followed the Einstein summation convention.
Using the Stokeslet, the velocity field of the fluid flow can be expressed as a single-layer boundary integral
\begin{equation}
    \bm{u}(\bm{x}) = \mathcal{S}[\bm{\xi}](\bm{x}), \quad \bm{x} \in \mathbb{R}^3,
\end{equation}
where $\bm{\xi} = \jump{\bm{\sigma}\bm{n}}$ represents an unknown single-layer density, introduced as an auxiliary variable. 
The auxiliary variable $\bm{\xi}$ has a clear physical interpretation as the exerted force, which can be directly used to enforce the force balance on $\Gamma$ in equation~\eqref{eqn:sf-force}.
The single-layer integral is continuous across the surface and its restriction on $\Gamma$ is well-defined, which we also denote by $\mathcal{S}[\bm{\xi}]$. Thus, we have $\bm{U} = \mathcal{S}[\bm{\xi}]$ on $\Gamma$.

Fixing the time $t = t^*$, the moving boundary problem can be regarded as a set of PDEs defined on a fixed surface.
Define $\mathfrak X: \Gamma\times [0,T]\to \mathbb R^3$ such that $\bm X(\bm \theta, t) = \mathfrak X(\bm X(\theta, t^*), t)$ and $\partial_t \mathfrak X = \partial_t \bm X \circ\bm X^{-1} = \bm{U}$ on $\Gamma$. Clearly, we have $\mathfrak X(\cdot, t^*) = \bm \Id(\cdot)$ where the identity mapping $\bm \Id$ is defined as $\bm \Id(\bm x) = \bm x$ for $\bm x \in \Gamma$. 
Due to the bending force in Equation~\eqref{eqn:sf-force}, which introduces the surface Laplacian of the mean curvature, the governing equations are high-order in space. 
To simplify, we treat the mean curvature $H$ as an unknown function and reformulate the equations in a mixed form. 
Additionally, we introduce auxiliary variables, including the mean curvature vector $\bm{H} = \frac{1}{2}\Delta_\Gamma \bm \Id$ and the normal bending force $g$. 
Then we reformulate the problem as finding $\bm{\xi}, \bm{U}, P, \bm{H}, H, g$ such that the following equations are satisfied,
\begin{subequations}\label{eqn:pde-mix-form}
\begin{align}
    - \mathcal{S}[\bm{\xi}] + \bm{U}  &= \bm{0}, \quad \text{ on }\Gamma,\\
    \bm{\xi}  - \nabla_{\Gamma}\cdot\bm{\sigma}_{\Gamma} +  g\bm{n} &= \bm{0},\quad \text{ on }\Gamma,\\
     \nabla_{\Gamma} \cdot \bm{U} &= \bm{0},  \quad \text{ on }\Gamma,\\
     \partial_t \mathfrak X   &=  \bm{U},  \quad \text{ on }\Gamma, \\
     2\bm{H} - \Delta_{\Gamma} \bm \Id  &= \bm{0},  \quad \text{ on }\Gamma, \\
     H - \bm{n}\cdot\bm{H} &= 0,  \quad \text{ on }\Gamma,  \\
     g - \alpha\paren{ \Delta_{\Gamma} H + 2\paren{H-c_0} \paren{H^2+Hc_0-K}} &= 0,  \quad \text{ on }\Gamma,
\end{align}
\end{subequations}
together with boundary conditions \eqref{eqn:bc-conormal}, \eqref{eqn:bc-normal}, and  \eqref{eqn:bc-moment}.

Let $\dual{\cdot}{\cdot}_\Gamma$ be the $L^2$ inner product defined on $\Gamma$.
We define the $L^2$ space and the Sobolev spaces on $\Gamma$ as
\begin{align}
    L^2\paren{\Gamma} = \left\{f: \dual{f}{f}_\Gamma < +\infty \right\},\quad
    H^1\paren{\Gamma} = \left\{f\in L^2\paren{\Gamma} , \grad_{\Gamma} f \in [L^2\paren{\Gamma}]^3\right\}.
\end{align}
For functions with Dirichlet boundary conditions on $\partial\Gamma$, we define
\begin{align}
    H_0^1\paren{\Gamma} &= \{f \in H^1(\Gamma): \at{f}{\partial\Gamma} = 0\},\\
    H_g^1 \paren{\Gamma} &= \left\{ f\in H^1\paren{\Gamma}: \at{f}{\partial\Gamma}=- \frac{\alpha_G}{\alpha} \kappa_n + c_0 \right\}.
\end{align}
In order to define the functional space for the single-layer density, we suppose that $\Gamma$ can be extended to a closed surface $S$. 
Let the $\dual{\cdot}{\cdot}_{S}$ be the dual pairing on $S$ defined by extending the $L^2$ inner product on $S$.
Define the trace space and its dual as
\begin{equation}
    H^{\frac{1}{2}}(S) = \{\at{f}{S}: f \in H^1(\mathbb{R}^3) \}, \quad H^{-\frac{1}{2}}(S) = (H^{\frac{1}{2}}(S))'.
\end{equation}
With the Sobolev space on the closed surface $S$, we can now define the fractional Sobolev space on the open surface $\Gamma$ by
\begin{equation}
    H^{\frac{1}{2}}(\Gamma) = \{\at{f}{\Gamma}: f\in H^{\frac{1}{2}}(S)\}, \quad \wt H^{-\frac{1}{2}}(\Gamma) = \{f\in H^{-\frac{1}{2}}(S): \text{supp }f \subset \overline \Gamma\}.
\end{equation}
Note that $\wt H^{\frac{1}{2}}(\Gamma)$ is the completion of $C_c^\infty(\Gamma)$ under the $H^{-\frac{1}{2}}(S)$ norm and it is the canonical dual of $H^{\frac{1}{2}}(\Gamma)$ with respect to dual pairing $\dual{\cdot}{\cdot}_S$ \cite{mclean2000strongly}.

The weak form of \eqref{eqn:pde-mix-form} is given as follows: for almost every $t\in[0,T]$, we find 
\begin{equation}
    \bm{\xi}  \in [\wt H^{-\frac{1}{2}}(\Gamma)]^3,\quad\bm{H}\in [L^2(\Gamma)]^3 , \quad
    \bm{U}, \mathfrak X  \in [H^1(\Gamma)]^3, \quad
    P, g\in L^2(\Gamma), \quad
    H \in H_g^1  (\Gamma).
\end{equation}
such that
\begin{subequations}\label{eqn:weak-form-3d}
\allowdisplaybreaks
\begin{align}
    &\dual{\bm{\vph}}{ -\mathcal{S}[\bm{\xi}] + \bm{U}}_{\Gamma}  = 0, \quad \forall \bm{\vph} \in [\wt H^{-\frac{1}{2}}(\Gamma)]^3,\\
    &\dual{\bm{\psi}}{\bm{\xi}}_{\Gamma} + 2\mu_\Gamma\dual{\bm{D}_{\Gamma}(\bm{\psi})}{\bm{D}_{\Gamma}(\bm{U})}_{\Gamma} - \dual{\nabla_{\Gamma}\cdot
    \bm{\psi}}{P}_{\Gamma} + \dual{\bm{\psi}}{g\bm{n}}_{\Gamma} 
    \\
    &\nonumber\, = -\dual{\bm{\psi}}{\paren{\alpha\paren{H-c_0}^2 + \alpha_G K + \gamma \kappa_g } \bm{\nu}}_{\partial\Gamma}, \quad \forall \bm{\psi}  \in [H^1(\Gamma)]^3,\\
    &-\dual{Q}{\nabla_{\Gamma}\cdot \bm{U}}_{\Gamma}  = 0,\quad \forall Q\in L^2\paren{\Gamma},\\
    &\dual{\bm{\omega}}{\partial_t \mathfrak X - \bm{U}}_{\Gamma}  =0, \quad \forall \bm{\omega}\in [L^2\paren{\Gamma}]^3,\\
   & 2\dual{\bm{\zeta}}{ \bm{H}}_{\Gamma} +  \dual{\nabla_{\Gamma}\bm{\zeta}}{\nabla_{\Gamma}\bm \Id }_{\Gamma}  = \dual{\bm{\zeta}}{ \bm{\nu}}_{\partial\Gamma}, \quad \forall \bm{\zeta}\in[H^1\paren{\Gamma}]^3, \label{eqn:zeta_h}\\
    &\dual{\eta}{H - \bm{n}\cdot\bm{H}}_{\Gamma}  = 0,\quad \forall \eta\in H_0^1\paren{\Gamma},\\
    &\dual{\chi}{g}_{\Gamma} + \alpha \dual{\nabla_{\Gamma}\chi}{\nabla_{\Gamma}H}_{\Gamma} - 2\alpha\dual{\chi}{\paren{H-c_0} \paren{H^2+Hc_0-K}}_{\Gamma} \\
   &\nonumber\, = \dual{\chi}{  \alpha_G \tau_g' - \gamma \kappa_n}_{\partial\Gamma},  \quad \forall \chi\in H^1\paren{\Gamma}.
\end{align}
\end{subequations}
It is important to note that, since the surface is open, the single-layer density $ \bm{\xi} $ is generally expected to belong to $\wt H^{-\frac{1}{2}}\paren{\Gamma}$, but not to $ L^2\paren{\Gamma} $ even when $\bm U \in H^1\paren{\Gamma}$. \cite{Stephan1987,Hayashi1977,Kirvalidze1997}. 

\subsection{Axisymmetric formulations}
While the full 3D formulation is applicable to arbitrary geometries, it is computationally intensive. To reduce computational complexity and enable more efficient numerical simulation, we only consider the axisymmetric case. 
This simplification assumes that the membrane surface is rotationally symmetric about the $ z $-axis and that the surrounding fluid flow is axisymmetric, with no azimuthal velocity component. As a result, the fluid velocity in the $ r$- and $z$-directions, along with the pressure, are taken to be independent of the azimuthal angle $\varphi$.

Let $\bm e_r, \bm e_\vph, \bm z$ be standard cylindrical basis vectors.
Let the surface $\Gamma$ be generated by the rotation of a planar curve $\mathcal{C}: (X^r(s), X^z(s))^{\top}$ around the $z$-axis, where $s\in\mathbb{T} = (0,L)$ denotes the arc-length parameter along the curve. 
The resulting axisymmetric surface $\Gamma$ is  parameterized as
\begin{equation}
\Gamma:\bm X (s,\varphi) = X^r(s)\bm{e}_r + X^z(s)\bm{e}_z 
, \quad s\in\mathbb T, \quad  \varphi \in\mathbb{S}^1 = \mathbb{R}\setminus (2\pi\mathbb{Z}).
\end{equation}
Since the surface is open, the 
generating curve $\mathcal{C}$ may terminate away from the  $z$-axis. We denote the endpoint lying on the $z$-axis by $\partial_0\mathcal{C}$, and the open edge by $\partial_1\mathcal{C}$. 
The mean and Gaussian curvatures are given by, respectively,
\begin{align}\label{eqn:axs-H-K}
    H = \frac{1}{2}\paren{X^r_s X^z_{ss} - X^z_s X^r_{ss} + \frac{X^z_s}{X^r}}, \quad
    K = \frac{X^z_s}{X^r} \paren{X^r_s X^z_{ss} - X^z_s X^r_{ss}}.
\end{align}
Note that the current configuration allows the surface to have up to two holes, which correspond to cases where both endpoints of the generating curve $\mathcal{C}$ do not intersect the $z$-axis. For such cases, we provide the expressions for the co-normal vector and curvatures on the boundary $\partial\Gamma$.
The co-normal vector to $\partial\Gamma$ is given by
\begin{align}
    \bm{\nu}  = 
    \begin{cases}
        X^r_s(L) \bm{e}_r + X^z_s(L) \bm{e}_z, & s = L,\\
        - X^r_s(0) \bm{e}_r - X^z_s(0) \bm{e}_z, & s = 0.
    \end{cases}
\end{align}
The normal and geodesic curvatures of $\partial\Gamma$ are given by, respectively,
\begin{align}
    \kappa_n =
    \begin{cases}
         \at{\frac{X^z_s}{X^r}}{s=L}, & s = L,\\
        -\at{\frac{X^z_s}{X^r}}{s=0} , & s = 0.
    \end{cases} ,\quad
    \kappa_g =
    \begin{cases}
         \at{\frac{X^r_s}{X^r}}{s=L} , & s= L,\\
         -\at{\frac{X^r_s}{X^r}}{s=0}, & s = 0.
    \end{cases}
\end{align}
Since the surface is axisymmetric, the torsion vanishes, i.e., $\tau_g = 0$.

For axisymmetric scalar and vector-valued functions $f$ and $\bm{f}$ defined on the surface $\Gamma$, the surface gradient and surface divergence simplify to the following forms
\begin{align}
    \nabla_{\Gamma} f =  f_s (X^r_s \bm{e}_r + X^z_s \bm{e}_z), \label{eqn:sur-grad0}\quad 
    \nabla_{\Gamma}\cdot\bm{f} = f^r_s X^r_s + f^z_s X^z_s + \frac{f^r}{X^r}.
\end{align}
Furthermore, the surface gradient and deformation tensor of $\bm{f}$ are given by
\begin{align}
    \nabla_{\Gamma}\bm{f} &= (f^r_s\bm{e}_r + f^z_s\bm{e}_z)\otimes (X^r_s \bm{e}_r +  X^z_s \bm{e}_z) + \frac{f^r}{X^r} \bm{e}_{\vph}\otimes\bm{e}_{\vph},\\
    \bm{D}_{\Gamma}\paren{\bm{f}} &= \paren{\bm X _s\cdot\bm{f}_s} (X^r_s \bm{e}_r +  X^z_s \bm{e}_z) \otimes (X^r_s \bm{e}_r +  X^z_s \bm{e}_z) + \frac{f^r}{X^r} \bm{e}_{\vph}\otimes\bm{e}_{\vph}.
\end{align}

The single-layer integral $\mathcal{S}[\bm{\xi}]$ on the axisymmetric surface $\Gamma$ reduces to
\begin{equation}
     \mathcal{S}[\bm{\xi}](\bm X \paren{s, \vph}) = - \frac{1}{\mu}\int_{\mathbb T}\int_{\mathbb S^1} G\paren{\bm X \paren{s,\vph} - \bm X \paren{s',\vph'} } \bm{\xi}\paren{\bm X \paren{s',\vph'}} X^r(s')  \,d\vph' \,ds',
\end{equation}
and admits an analytic integration in the $\vph$-direction, yielding the axisymmetric single-layer representation for the surface velocity $\bm U$,
\begin{align}\label{eqn:axs-sin-int}
    \bm{U}(s) =
    \begin{pmatrix}
        (\mathcal{S}[\bm{\xi}])^r \\
        (\mathcal{S}[\bm{\xi}])^z
    \end{pmatrix}(s) 
    = \frac{1}{8\pi\mu} \int_{\mathbb{T}}  X^r\paren{s'} 
    \bm{S}(s,s')
    \bm{\xi}(s')
     \,ds' \eqdef  \mathcal{S}[\bm{\xi}](s),\quad s\in\mathbb T,
\end{align}
where both components $(\mathcal{S}[\bm{\xi}])^r$ and $(\mathcal{S}[\bm{\xi}])^z$ are independent of $\vph$. 
The kernel $\bm{S}\paren{s, s'}$ is given by
\begin{align}\label{eqn:axs-ker}
    \bm{S}(s,s') = 
    \begin{pmatrix}
      I_{11} + (r^2+(r')^2)I_{31} - r r'(I_{30}+I_{32}) & \Delta z(rI_{30} - r'I_{31})\\
       \Delta z (rI_{31} - r'I_{30}) & I_{10}+(\Delta z)^2I_{30}
    \end{pmatrix},
\end{align}
with $r = X^r(s), r' = X^r(s'), \Delta z = X^z(s)-X^z(s')$. 
Here, the function $I_{mn} = I_{mn}(r, r', \Delta z)$ is given by
\begin{align}\label{eqn:Imn}
    I_{mn}(r_1, r_2, \hat{z}) = \frac{2}{((r_1-r_2)^2 + \hat{z}^2)^{\frac{m}{2}}} I_{mn}^{(1)}(k_1) + \frac{2}{((r_1+r_2)^2 + \hat{z}^2)^{\frac{m}{2}}}  I_{mn}^{(2)}(k_2),
\end{align}
where $I_{mn}^{(1)}(k)$ and $I_{mn}^{(2)}(k)$ are related to complete elliptic integrals,
\begin{align}
&I_{mn}^{(1)}(k_1) = \int_0^{\frac{\pi}{2}}  \frac{(1-2\sin^2\theta')^n}{(1 + k_1^2 \sin^2\theta')^{\frac{m}{2}}}  \,d\theta', \quad k_1^2 = \frac{4r_1r_2}{(r_1-r_2)^2 + \hat{z}^2},\\
&I_{mn}^{(2)}(k_2) = \int_{0}^{\frac{\pi}{2}}  \frac{(2\sin^2\theta'-1)^n}{(1 - k_2^2 \sin^2\theta')^{\frac{m}{2}}}  \,d\theta',\quad k_2^2 = \frac{4 r_1r_2}{(r_1+r_2)^2 + \hat{z}^2} .
\end{align}
For example, $I_{10}^{(2)} (k) =  K(k)$ and $I_{11}^{(2)} (k) =  \frac{2}{k^2}(K(k)-E(k)) - K(k)$ where $K(k)$ and $E(k)$ are complete elliptic integrals of the first and second kind, respectively.

The $L^2$ inner-product $\dual{\cdot}{\cdot}_\Gamma$ naturally defines the weighted $L^2$ inner-product on the generating curve $\mathcal{C}$,
\begin{align}
    \dual{f}{g}_\Gamma = 2\pi \int_{\mathbb{T}} X^r\paren{s} f\paren{s}g\paren{s} \,ds \eqdef 2\pi\dual{f}{g}_{\mathcal{C}}.
\end{align} 
The weighted $L^2$ and Sobolev spaces on $\mathcal{C}$ are defined as
\begin{align}
    L^2\paren{\mathcal{C}} = \left\{f:\dual{f}{f}_{\mathcal{C}} < +\infty\right\}, \quad
    H^1\paren{\mathcal{C}} = \left\{f\in L^2\paren{\mathcal{C}}, \partial_s f \in L^2\paren{\mathcal{C}}\right\}.
\end{align}
Similarly, we define
\begin{align}
    H_0^1 \paren{\mathcal{C}} &= \left\{ f\in H^1\paren{\mathcal{C}}: \at{f}{\partial_1\mathcal{C}} =0 \right\}, \\
    H_g^1 \paren{\mathcal{C}} &= \left\{ f\in H^1\paren{\mathcal{C}}: \at{f}{\partial_1\mathcal{C}} = - \frac{\alpha_G}{\alpha} \kappa_n + c_0 \right\}.
\end{align}
The following functional space is defined for vector-valued axisymmetric functions
\begin{align}\label{eqn:axis-space}
    V_{\partial}\paren{\mathcal{C}} = \left\{\bm{f} \in [H^1(\mathcal{C})]^2 : \at{f^r}{\partial_0 \mc C} = 0\right\}.
\end{align}
The negative-order Sobolev space $[\wt H^{-\frac{1}{2}}(\mc C)]^2$ is defined as
\begin{equation}
    [\wt H^{-\frac{1}{2}}(\mc C)]^2 = \{(f^r, f^z): f^r\bm e_r + f^z \bm e_z \in \wt H^{-\frac{1}{2}}(\Gamma)\}.
\end{equation}
We now formulate the weak form under the axisymmetric setting as follows: for almost every $t\in[0,T]$, we find
\begin{align}
    \bm{\xi} \in [\wt H^{-\frac{1}{2}}\paren{\mathcal{C}}]^2 ,\quad \bm{H} \in [L^2\paren{\mathcal{C}}]^2 , \quad
    \bm{U},\mathfrak X\in V_{\partial}\paren{\mathcal{C}}, \quad
    P, g\in L^2\paren{\mathcal{C}}, \quad
    H \in H_g^1  \paren{\mathcal{C}},
\end{align}
such that
\begin{subequations}\label{eqn:weak-axisym}
\allowdisplaybreaks
\begin{align}
    &\dual{\bm{\vph}}{ -\mathcal{S}[\bm{\xi}] + \bm{U}}_{\mathcal{C}} = 0,\quad\forall \bm{\vph}\in [\wt H^{-\frac{1}{2}}(\mc C)]^2,\\
    &\dual{\bm{\psi}}{\bm{\xi}}_{\mathcal{C}} +2\mu_\Gamma \dual{\bm X _s\cdot\bm{\psi}_s}{\bm X _s\cdot\bm{U}_s}_{\mathcal{C}} +2\mu_\Gamma\dual{\frac{\psi^r}{X^r}}{\frac{U^r}{X^r}}_{\mathcal{C}} \\
    &\nonumber\, - \dual{\bm X _s\cdot\bm{\psi}_s + \frac{\psi^r}{X^r}}{P}_{\mathcal{C}} + \dual{\bm{\psi}}{g\bm{n}}_{\mathcal{C}} \\
    &\nonumber\, = -\dual{\bm{\psi}}{\paren{\alpha\paren{H-c_0}^2 + \alpha_G K + \gamma \kappa_g } \bm{\nu}}_{\partial\mathcal{C}}, \quad \forall\bm{\psi}\in V_{\partial}\paren{\mathcal{C}},\\
    &-\dual{Q}{\bm X _s\cdot\bm{U}_s + \frac{U^r}{X^r}}_{\mathcal{C}} = 0, \quad\forall Q\in L^2\paren{\mathcal{C}},\\
    &\dual{\bm{\omega}}{\partial_t \mathfrak X  - \bm{U}}_{\mathcal{C}} =0, \quad \forall\bm{\omega}\in [L^2\paren{\mathcal{C}}]^2,\\
    &2\dual{\bm{\zeta}}{ \bm{H}}_{\mathcal{C}} + \dual{\bm{\zeta}_s}{\bm X _s}_{\mathcal{C}} + \dual{\frac{\zeta^r}{X^r}}{\frac{X^r}{X^r}}_{\mathcal{C}} = \dual{\bm{\zeta}}{ \bm{\nu}}_{\partial\mathcal{C}}, \quad\forall\bm{\zeta}\in V_{\partial}\paren{\mathcal{C}},\\
    &\dual{\eta}{H - \bm{n}\cdot\bm{H}}_{\mathcal{C}}= 0, \quad \eta\in H_0^1\paren{\gamma},\\
    &\dual{\chi}{g}_{\mathcal{C}} + \alpha \dual{\chi_s}{H_s}_{\mathcal{C}} - 2\alpha\dual{\chi}{\paren{H-c_0} \paren{H^2+Hc_0-K}}_{\mathcal{C}} \\
    &\nonumber\,= \dual{\chi}{\alpha_G \tau_g' - \gamma \kappa_n}_{\partial\mathcal{C}}, \quad \forall\chi \in H^1\paren{\mathcal{C}}.
\end{align}
\end{subequations}

\section{A hybrid BEM-FEM method}\label{sec:4}
\subsection{Finite element discretization}

{\color{black}
We next construct a fully discrete scheme for the weak formulation~\eqref{eqn:weak-axisym}. We discretize time using a semi-implicit backward Euler scheme, and we discretize space using a Galerkin finite element method on the generating curve.

At each time level $t^{n}$, the membrane configuration is represented by the generating curve $\mc C^n$, given by the numerical approximation $\bm X^n(\cdot)$. The only time derivative in the system appears in the curve evolution $\partial_t \bm X = \bm U$, which we discretize using the finite  difference $\paren{\bm X^{n+1} - \bm X^{n}}/\Delta t.$
For simplicity, we adopt uniform time stepping $t^n = n \Delta t$ for $n = 0,1,2,\dots,N_t$. All unknown quantities in~\eqref{eqn:weak-axisym} are then understood as functions defined on $\mc C^n$. To avoid remeshing on a moving geometry, we pull these quantities back to a fixed reference interval $\mathbb I = (0,1)$ using the current parametrization $\bm X^n$. For any surface quantity $f : \mc C^n \to \mathbb{R}$ we denote by $(\bm X^{n,*} f)(\alpha) := f(\bm X^{n}(\alpha))$ its pullback to $\mathbb I$.}
We then approximate these pulled-back unknowns in finite element spaces on $\mathbb I$. The reference interval $\mathbb I$ is partitioned into $N$ cells with nodes $0 = \alpha_0 < \alpha_1 < \cdots < \alpha_N = 1.$
On this mesh, we define a one-dimensional continuous piecewise polynomial finite element space as
\begin{equation}
\begin{aligned}
    V^N_k(\mathbb{I}) &= \left\{
        v \in C(\overline{\mathbb{I}}) :
        \left.v\right|_{[\alpha_{i-1}, \alpha_i]} \in \mathbb{P}_k([\alpha_{i-1}, \alpha_i]),
        \; 1 \le i \le N
    \right\} \\
    &= \operatorname{span} \left\{
        \phi_k^1(x), \phi_k^2(x), \dots, \phi_k^{M_k}(x)
    \right\},
\end{aligned}
\end{equation}
where $\mathbb{P}_k$ is the space of polynomials of degree at most $k$. We define $V_{\partial}(\mathbb{I})$, $H_0^1(\mathbb I)$, and $H_g^1(\mathbb I)$ analogously to $V_{\partial}(\mc C)$, $H_0^1(\mc C)$, and $H_g^1(\mc C)$, respectively.

Let $\bm X^n(\cdot) \in [V^N_2(\mathbb{I})]^2 \cap V_{\partial}(\mathbb{I})$ denote the numerical approximation to $\bm X(\cdot,t^n)$, and let $f^{n}(\cdot)\in V_k^N(\mathbb I)$ denote the numerical approximation to the pullback $(\bm X^{n,*} f)(\cdot) = f(\bm X^{n}(\cdot)).$
To ensure the inf--sup stability condition for the surface Stokes equations, we employ a $P_2$--$P_1$ Taylor--Hood element: quadratic elements for the surface velocity $\bm U$ and linear elements for the surface tension $P$. For the remaining quantities, we use quadratic finite element spaces for their approximations.

Let $ \bm{\xi}^n $, $ \bm{U}^n $, $ P^n $, $ \bm{H}^n $, $ H^n $, $ g^n $ denote the numerical approximations to the pullbacks of $ \bm{\xi} $, $ \bm{U} $, $ P $, $ \bm{H} $, $ H $, $ g $ with respect to $\bm X^n$, respectively.
The initial condition $\bm X^0(\cdot)$ is obtained by the interpolation of $\bm X(\cdot, 0)$.
The hybrid BEM-FEM method is formulated as: for $n\geq 0$, we find
\begin{align}
    &\bm{H}^{n+1}, \bm{\xi}^{n+1}\in [V^N_2(\mathbb{I})]^2 , \quad \bm{U}^{n+1}, \bm X ^{n+1} \in [V^N_2(\mathbb{I})]^2 \cap V_{\partial}(\mathbb{I}), \\
    &P^{n+1} \in V^N_1(\mathbb{I}), \quad  g^{n+1} \in V^N_2(\mathbb{I}), \quad
    H^{n+1} \in V^N_2(\mathbb{I}) \cap H^1_g(\mathbb{I}),
\end{align}
such that
\begin{subequations}\label{eqn:scheme}
\allowdisplaybreaks
\begin{align}
    & - \dual{\bm{\vph}}{ \mathcal{S}^{h,n}[\bm{\xi}^{n+1}]}_{\mathcal{C}^n}^{h,\ast} + \dual{\bm{\vph}}{ \bm{U}^{n+1}}_{\mathcal{C}^n}^h = 0,\forall \bm \vph \in [V^N_2\paren{\mathbb{I}}]^2 ,\\
    &\dual{\bm{\psi}}{\bm{\xi}^{n+1}}_{\mathcal{C}^n}^h +2\mu_\Gamma \dual{\bm X _s^n\cdot\bm{\psi}_s}{\bm X _s^n\cdot\bm{U}_s^{n+1}}_{\mathcal{C}^n}^h +2\mu_\Gamma\dual{\frac{\psi^r}{X^{r,n}}}{\frac{U^{r,n+1}}{X^{r,n}}}_{\mathcal{C}^n}^h \\
    &\nonumber\, - \dual{\bm X _s^n\cdot\bm{\psi}_s + \frac{\psi^r}{X^{r,n}}}{P^{n+1}}_{\mathcal{C}^n}^h +\dual{\bm{\psi}}{g^{n+1}\bm{n}^n}_{\mathcal{C}^n}^h \\
    &\nonumber\, = -\dual{\bm{\psi}}{\paren{\alpha\paren{H^n-c_0}^2 + \alpha_G K^n + \gamma \kappa_g^n } \bm{\nu}^n}_{\partial\mathcal{C}^n}, \forall\bm{\psi}\in [V^N_2\paren{\mathbb{I}}]^2 \cap V_{\partial}\paren{\mathbb{I}},\\
    &-\dual{Q}{\bm X _s^n\cdot\bm{U}_s^{n+1} + \frac{U^{r,n+1}}{X^{r,n}}}_{\mathcal{C}^n}^h = 0, \forall Q\in V^N_1\paren{\mathbb{I}},\\
    & \dual{\bm{\omega}}{\frac{1}{\tau}\paren{\bm X ^{n+1} - \bm X ^n} - \bm{U}^{n+1}}_{\mathcal{C}^n}^h =0,  \forall\bm{\omega}\in[V^N_2\paren{\mathbb{I}}]^2,\\
    &{\color{black}  2\dual{\bm{\zeta}}{ \bm{H}^{n+1}}_{\mathcal{C}^n}^h + \dual{\bm{\zeta}_s}{\bm X _s^{n+1}}_{\mathcal{C}^n}^h + \dual{\frac{\zeta^r}{X^{r,n}}}{\frac{X^{r,n+1}}{X^{r,n}}}_{\mathcal{C}^n}^h }\\
    &{\color{black} \nonumber\, = \dual{\bm{\zeta}}{ \bm{\nu}^n}_{\partial\mathcal{C}^n}, \forall\bm{\zeta}\in [V^N_2\paren{\mathbb{I}}]^2 \cap V_{\partial}\paren{\mathbb{I}},}\\
    & \dual{\eta}{H^{n+1} - \bm{n}^n\cdot\bm{H}^{n+1}}_{\mathcal{C}^n}^h= 0,  \forall\eta \in V^N_2\paren{\mathbb{I}} \cap H^1_0\paren{\mathbb{I}}, \\
    &\dual{\chi}{g^{n+1}}_{\mathcal{C}^n}^h + \alpha \dual{\chi_s}{H_s^{n+1}}_{\mathcal{C}^n}^h = 2\alpha\dual{\chi}{\paren{H^n-c_0} \paren{\paren{H^n}^2+H^n c_0-K^n}}_{\mathcal{C}^n}^h  \\
    &\nonumber\,+ \dual{\chi}{ - \gamma \kappa_n^n}_{\partial\mathcal{C}^n}, \forall\chi \in V^N_2\paren{\mathbb{I}}.
\end{align}
\end{subequations}
Here, the discrete inner-product $\dual{\cdot}{\cdot}_{\mathcal{C}^n}^h$ is an approximation to the inner product $\dual{\cdot}{\cdot}_{\mathcal{C}^n}$ by Gauss-Legendre quadrature, and the double-integral $\dual{\cdot}{ \mathcal{S}^{h,n}[\cdot]}_{\mathcal{C}^n}^{h,\ast} $ is computed based on the Alpert’s quadrature \cite{Alpert1999}, which is well-suited for handling logarithmic-type weakly singular integrals.

Due to the presence of high-order spatial derivatives, any explicit scheme suffers from numerical stiffness, imposing a severe constraint on the time step. While fully implicit schemes are generally more stable and can sometimes be unconditionally stable, the nonlinear nature of the equation makes such schemes challenging, as they result in nonlinear algebraic systems that are difficult and time-consuming to solve. The current numerical scheme \eqref{eqn:scheme} uses a semi-implicit time discretization designed to balance stability and efficiency. In this scheme, only the geometry and boundary conditions are treated explicitly, while terms with high-order spatial derivatives are treated implicitly, thus removing the severe time-step constraint and yielding a linear system.


The resulting coefficient matrix is mostly sparse, except for the block corresponding to the single-layer integral, which remains dense. However, due to the axisymmetric reduction, the problem is effectively one-dimensional, keeping the size of the linear system small. This enables efficient solution using direct methods, such as \texttt{SparseLU} \cite{eigenweb}.

\subsection{Singularity capturing mesh refinement}
The solution to the open membrane problem typically exhibits singular behavior at the open edge.
{\color{black}The single-layer density diverges as $ |\bm \xi| \sim d^{-\frac{1}{2}} $,} where $ d $ is the distance to the edge. 
Due to fluid-membrane coupling, membrane tension, velocity, and the curve itself can also exhibit certain singularity, typically weaker than $\bm \xi$.
Such singularities can easily deteriorate the convergence of the finite element method if a uniform computational mesh is employed on the curve. 
Here, we present a direct method for introducing local mesh refinement to capture the singular behavior at the edge by only modifying the initial parametrization of the generating curve.

Let $ \bm X^0(\alpha), \alpha \in \mathbb I $ represent the initial parameterization proportional to the arc-length $s$, i.e., $|\partial_\alpha \bm X^0| (\alpha) \equiv C$ where $C>0$ is a constant.
Assuming the curve is sufficiently smooth, a uniform partition in $ \alpha $ results in almost uniform spacing of points along the curve. By carefully selecting the reparametrization, a uniform mesh in the new parameter space can produce locally refined points along the curve, which is particularly useful for capturing singular behavior in the solution.
Let $\Phi(\cdot):\mathbb I \to \mathbb I$ be a locally refined reparametrization function that satisfies the following properties:
\begin{enumerate}
    \item The function $\Phi$ maps the two endpoints $0$ and $1$ to themselves, i.e., $ \Phi(0) = 0 $, $ \Phi(1) = 1 $;
    \item The reparametrization is non-singular, i.e., $ \Phi' (\eta) > 0 $ for $ \eta \in (0,1) $;
    \item The reparametrization $\Phi$ has a local refinement property near the open edge, i.e., $ \Phi' (\eta) $ is close to $ 0 $ for all $ \eta $ such that $ \bm X^0\paren{\Phi(\eta)} \in \partial_0 \mathcal{C} $.
\end{enumerate}
If $\bm X^0(0)$ represents the closed end and $ \bm X^0(1) $ represents the open end, we consider the following function,
\begin{align}
    \Phi(\eta) =\cos\left(\frac{\pi}{2}(1-\eta)\right).
\end{align}
It can be verified that $ \Phi $ satisfies the three properties. 
let $N$ be the total number of mesh cells. {\color{black}By uniformly partitioning} $\eta$ with $\Delta \eta = 1/N$, one obtains the following mesh
\begin{equation}
    0 =\alpha_0 < \alpha_1 < \cdots < \alpha_N = 1, \quad \alpha_i = \Phi(i\Delta \eta), i = 0,1,\cdots, N.
\end{equation}
This is a graded mesh with points increasingly clustered near the open endpoint.
In particular, let $h_i=\alpha_{i+1}-\alpha_i$. The mesh spacing is approximately $h_i=\mathcal{O}(N^{-2})$ near the open endpoint and $h_i=\mathcal{O}(N^{-1})$ elsewhere.
A graded mesh of this type admits the following approximation property: for any solution $\bm\psi \in [\wt H^{-\frac{1}{2}}(\mc C)]^2$, there exists an approximation $\bm\psi_h = \wt{\bm\psi}_h \circ (\bm X^0)^{-1}$ with $\wt{\bm\psi}_h \in [V_k^N(\mathbb I)]^2$, and some $\delta > 0$ such that
\begin{equation}
\norm{\bm\psi - \bm \psi_h}_{\wt H^{-\frac{1}{2}}(\mc C)} \le C N^{-(1-\delta)},
\end{equation}
where the constant $C$ is independent of $N$~\cite[Lemma 3.1]{VonPetersdorff1990}.

Note that this method results in poor parametrization near the endpoints, where $\Psi’(\eta) = 0$, which may lead to numerical instability due to the explicit handling of geometric quantities such as $\kappa_n$ and $\kappa_g$ in the boundary conditions.
A simple remedy is to use a regularized version
\begin{align}
    \wt\Phi(\eta) = \paren{1-\varepsilon} \Phi(\eta) + \varepsilon \eta,
\end{align}
where $ \varepsilon > 0 $ is a small regularization parameter, typically chosen as $ \varepsilon = 10^{-3} $. This regularization is not always necessary and is introduced primarily for numerical convenience. Since the geometric quantities at the endpoints are intrinsically independent of the parametrization, they can alternatively be computed using a regular local parametrization of the generating curve. In such cases, one may simply set $ \varepsilon = 0 $.

\section{Numerical experiments}\label{sec:5}
In this section, we present numerical experiments to validate our numerical method for simulating an open membrane interacting with a surrounding fluid flow.
To ensure generality and facilitate comparison with related studies, we introduce the nondimensional form of the governing equations. Let $L$ and $T$ represent the characteristic length and time scales, respectively. We define the following nondimensional variables
\begin{equation}
\begin{aligned}
    \widetilde{\bm{x}} = \frac{\bm{x}}{L}, \quad \widetilde{t} = \frac{t}{T}, \quad \widetilde{\bm{u}} = \frac{T\bm{u}}{L}, \quad \widetilde{\bm{U}} = \frac{T\bm{U}}{L}, \quad \widetilde{p} = \frac{T p}{\mu}, \quad \widetilde{P} = \frac{T P}{\mu_{\Gamma}}.
\end{aligned}
\end{equation}
Setting the time scale as $T = L^3 \mu / \alpha$ and omitting the tildes for clarity, the nondimensional equations are given by
\begin{subequations}\label{eqn:nondim-pde}
\allowdisplaybreaks
\begin{align}
     -\Delta \bm u + \grad p &= \bm{0},\quad \text{ in } \Omega,\\
    -\nabla\cdot\bm{u} &= 0,\quad \text{ in } \Omega,\\
    \jump{\bm{\sigma}\bm{n}} - \beta \nabla_{\Gamma}\cdot\bm{\sigma}_{\Gamma} + \paren{ \Delta_{\Gamma} H + 2\paren{H-H_0} \paren{H^2+H_0 H-K}}\bm{n} &= \bm{0},\quad \text{ on }\Gamma,\\
     -\nabla_{\Gamma} \cdot \bm{U} &= \bm{0},  \quad \text{ on }\Gamma,\\
     \bm{u}   &= \bm{U},  \quad \text{ on }\Gamma,\\
     \beta\bm{\sigma}_{\Gamma}\bm{\nu} +  \paren{ \paren{H- H_0}^2 + \gamma_g K + \gamma_l\kappa_g } \bm{\nu}  &= \bm{0},\quad \text{ on }\partial\Gamma,\\
      \nabla_{\Gamma}H\cdot\bm{\nu} - \gamma_g \tau_g' + \gamma_l \kappa_n &= 0, \quad \text{ on }\partial\Gamma,\\
     H - H_0 + \gamma_g \kappa_n &= 0,\quad \text{ on }\partial\Gamma,\label{eqn:nondim-bc-3}
\end{align}
\end{subequations}
together with the surface evolution equation $\partial_t \bm X\circ\bm X^{-1} = \bm U$ and the nondimensional parameters
\begin{align}
    \beta = \frac{\mu_{\Gamma}}{L\mu}, \quad \gamma_g = \frac{\alpha_G}{\alpha},\quad \gamma_l = \frac{ \gamma L}{\alpha}, \quad H_0 = L c_0.
\end{align}

\subsection{Annular-shaped planar membrane}
We consider the case where the membrane is planar and has an annular shape,
\begin{equation}\label{eqn:planar-Gamma}
    \Gamma = \left\{(x,y,z) : z = 0, r = \sqrt{x^2 + y^2} \in \paren{R_i, R_o} \right\},\quad R_o > R_i > 0,
\end{equation}
where $R_i, R_o$ are the inner and outer radii of the membrane, respectively.
For such a configuration, our model reduces to a set of equations that are similar to the one in \cite{Jia2022}. Particularly, {\color{black}when the outer radius of the annular goes to infinity}, closed form solutions for the velocity field and the membrane tension can be constructed analytically, allowing us to validate our numerical method by comparing with the analytical solution. In this case, the bending energy of the membrane is ignored and the dynamics of the membrane is only driven by the line tension on the open edges. 

\begin{figure}[htbp]
    \centering
    \includegraphics[width=0.35\linewidth]{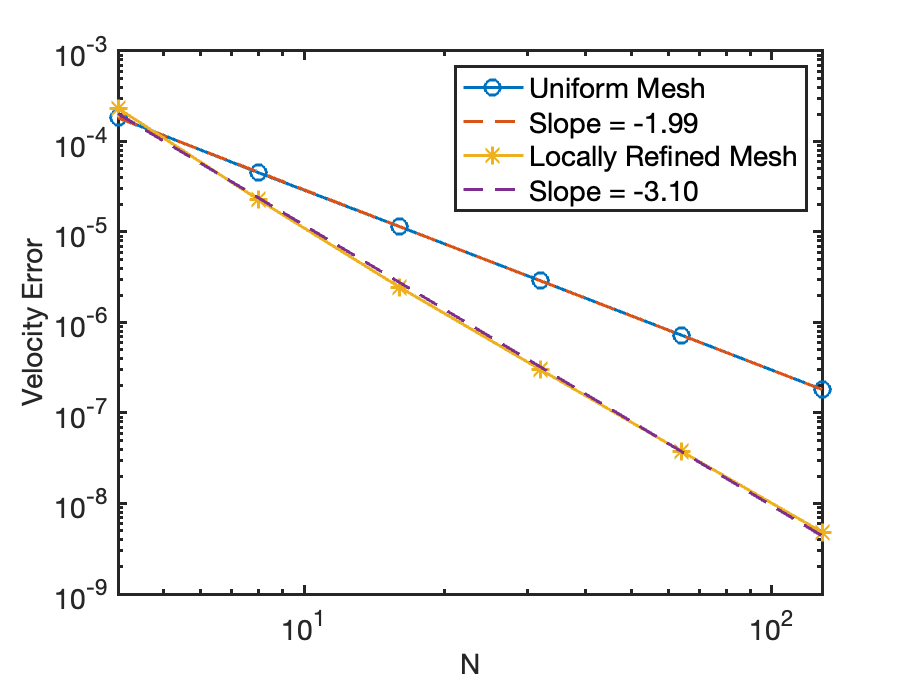} 
    \includegraphics[width=0.35\linewidth]{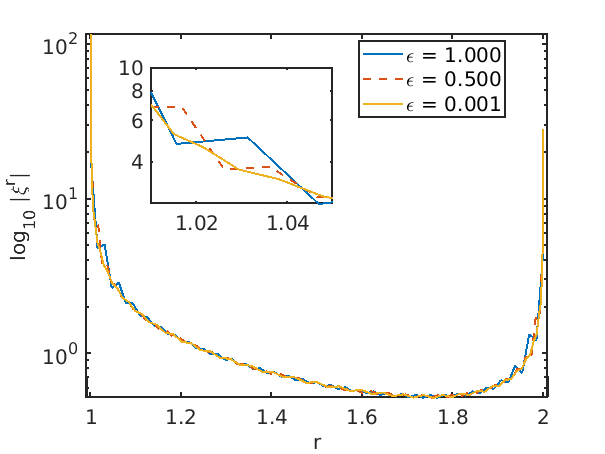} 
    \caption{Left: numerical error of the radial velocity for different mesh size using both uniform mesh and locally refined mesh; Right: radial single-layer density $\xi^r$ computed using different regularizing parameters $\varepsilon = 1, 0.5, 0.001$ with $N=32$.}
    \label{fig:eg1-1}
\end{figure}

We first perform a convergence test for the planar membrane case. Since the single-layer density is essentially unbounded, we evaluate only the error in the velocity.
The parameters used are chosen as $\beta = \gamma_l = 1$, and the membrane is defined with initial radii $ R_{i,0} = 1 $ and $ R_{o,0} = 2 $. Due to the inextensible condition, the planar membrane only has nonzero radial velocity, which is given by
$
    U^r(r) = \frac{F}{r}
$
where $F$ is an integration constant that can be explicitly computed when the outer radius of the membrane is infinite. For a finite membrane, the constant $F$ can be computed as follows. We numerically solve the PDE for one time step to obtain the numerical solution $U^r_h$. On noting that $F \equiv r U^r(r)$, we estimate the constant $ F $ by
\begin{equation}
    F_h = \frac{\int_{\Gamma} r U^r_h  \,dA}{|\Gamma|} = \frac{\int_{R_i}^{R_o} r^2 U^r_h(r)\,dr}{\int_{R_i}^{R_o} r\,dr}.
\end{equation}
Next, the error is defined as
\begin{equation}
    \norm{e_{U,h}}_{L^2\paren{\Gamma}} = \paren{\int_{R_i}^{R_o}r \paren{U^r_h(r) - \frac{F_h}{r} }^2 \, dr }^{\frac{1}{2}},
\end{equation}
where all the integrals are evaluated using the Gauss quadrature.
We compare the numerical errors obtained using a uniform mesh and a locally refined mesh for various mesh sizes $ N = 4, 8, \ldots, 128 $. The results, shown in the left panel of Figure~\ref{fig:eg1-1}, demonstrate that the uniform mesh degrades the convergence of the P2 element, whereas the locally refined mesh yields optimal convergence. 
We also vary the parameter $ \varepsilon $ to assess the scheme’s ability to capture the singularity of the single-layer density $ \bm{\xi} $. In the planar configuration, $ \bm{\xi} $ has only a nonzero radial component $ \xi^r $. Using a fixed mesh size $ N = 32 $, we plot the absolute value of $ \xi^r $ in the right panel of Figure~\ref{fig:eg1-1}.
When the mesh is uniform, corresponding to $\varepsilon=1$, the result generates noticeable oscillation in the vicinity of the open edge where $\xi^r$ is expected to diverge. On the contrary, when local mesh refinement is applied, corresponding to $\varepsilon=0.001$, oscillation near the open edge is suppressed and the predicted value at $x=1$ is much larger compared with uniform mesh, suggesting the local mesh refinement strategy increases the capability of capturing the singularity behavior of the solution near the open edge.  

\begin{figure}[htbp]
    \centering
    \includegraphics[width=0.35\linewidth]{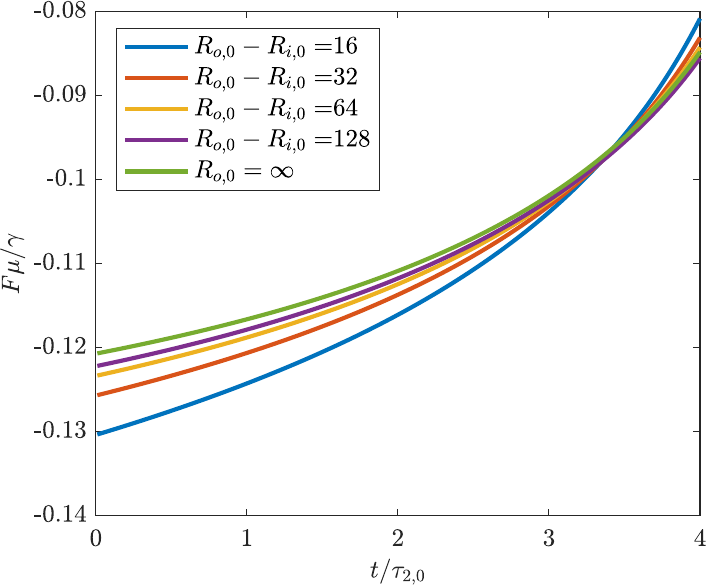}
    \caption{Numerical solution of $F$, normalized by $\gamma/\mu$, for different membrane widths $R_{o,0}-R_{i,0}=16,\,32,\,64,\,128$, compared with the analytical solution for $R_{o,0}=\infty$.}
    \label{fig:F-width}
\end{figure}

As the outer radius $R_o$ goes to infinity, an explicit expression for the value of the constant $F$ is available. Since we can only apply a finite number of $ R_o $ in the numerical method, we fix the initial inner radius as $R_{i,0}=1$ and successively increase the initial outer radius $ R_{o,0} $ to approximate the case where $ R_o = \infty $, and compare the numerical result of $ F(t) $ with the analytical result. We set the computation mesh as $ N = 32 $ with local mesh refinement for $\varepsilon=0.001$ and time step as $ \Delta t = 0.01 $ for the computation. The result is shown in Figure~\ref{fig:F-width}. The time is rescaled to $t/\tau_{2,0}$ where $\tau_{2,0} = \mu R_{i,0}^2/\gamma$ is the same as the one defined in \cite{Jia2022}. {\color{black}It can be observed that as the value $R_{o,0}-R_{i,0}$ increases, the estimated value of $F$ is closer to the limiting case where $R_{o,0}=\infty$}, demonstrating the convergence of the numerical solution. 



\begin{figure}[htbp]
    \centering
    \includegraphics[width=0.35\textwidth]{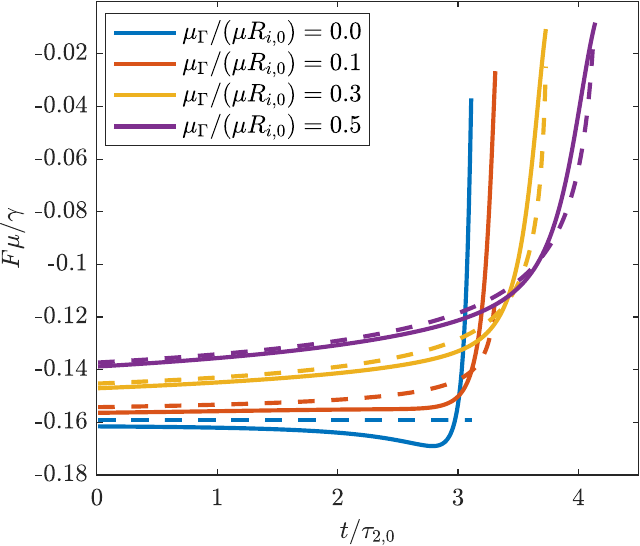}
    \includegraphics[width=0.35\textwidth]{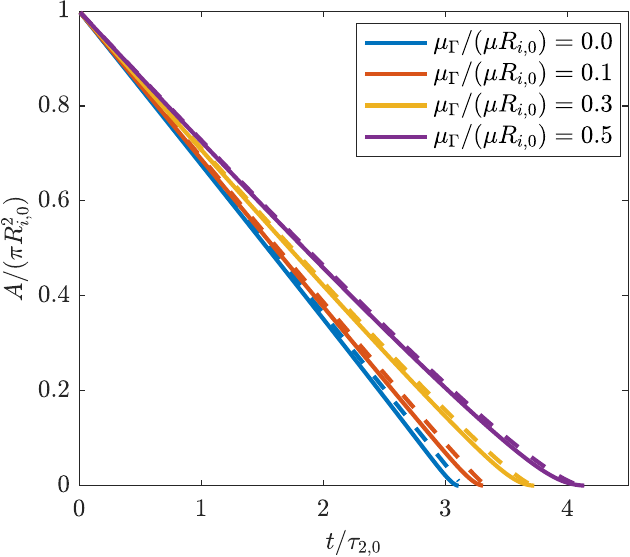}
    \caption{Comparison of numerical results (solid) and analytical solutions (dashed) for the value $F$, normalized by $\gamma/\mu$, and the dimensionless hole area $A/(\pi R_{i,0}^2)$ over time, for viscosity ratios $\mu_{\Gamma}/(\mu R_{i,0}) = 0,\,0.1,\,0.3,\,0.5$.}
    \label{fig:Fvst_Avst}
\end{figure}

{\color{black} We use the same computational setup to simulate hole-closing dynamics for several viscosity ratios $\mu_{\Gamma}/(\mu R_{i,0})$. The computed value $F$ and hole area $A$ are compared with the exact solutions obtained using the method of~\cite{Jia2022}, see Figure~\ref{fig:Fvst_Avst}. The numerical error increases as $\mu_{\Gamma}/(\mu R_{i,0})$ decreases. Moreover, as time progresses, the pore shrinks, strengthening the edge singularity (cf.~\cite{Jia2022}) and leading to larger numerical errors near closure. In particular, the case $\mu_{\Gamma}/(\mu R_{i,0})=0$ corresponds to a degenerate limit in which the pore dynamics reduce to a curve-shortening flow. The discrepancy between the numerical and exact results for $F$ around $t/\tau_{2,0} \approx 3.1$ reflects the increased difficulty in accurately resolving this near-closure regime. }

\subsection{Equilibrium Shapes}
Here we compute the equilibrium shapes of an open membrane and compare the results with those reported in \cite{Saitoh1998} and \cite{Tu2010}.  We simulate the relaxation of an inextensible elastic open membrane to equilibrium where the fluid flow ceases ($\bm{u} = 0$ and $\bm{U} = 0$) and the membrane shape $\Gamma$ is independent of time. Under the assumption of membrane inextensibility, the membrane area remains constant during the relaxation process, and is solely determined by its initial configuration. Consequently, the equilibrium shape depends on the initial shape and the following key parameters: $\gamma_g = \alpha_G/\alpha$, $\widetilde{\gamma} = \gamma/\alpha$, and $c_0$.

It is worth noting that our inextensibility assumption for the membrane differs from the extensible membrane model in \cite{Tu2010}, where a surface energy term associated with surface tension is included. In their formulation, surface tension is prescribed as a constant parameter, whereas in our approach, it emerges naturally as a Lagrange multiplier to enforce inextensibility. However, at equilibrium, both models have a constant surface tension across the membrane, enabling meaningful comparisons of equilibrium shapes.

\begin{figure}[htbp]
    \centering
    \includegraphics[width=\linewidth]{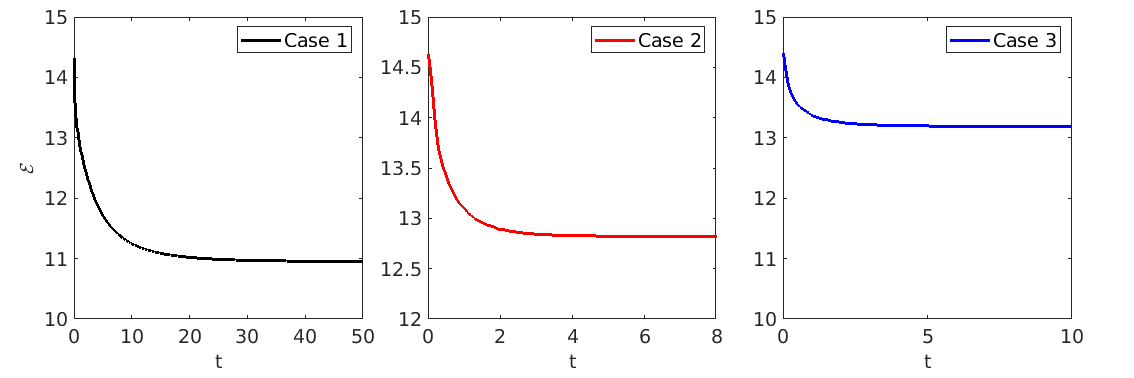}
    \caption{Time evolution of the total energy of the membrane for the three cases.}
    \label{fig:e2-energy}
\end{figure}

\begin{figure}[htbp]
    \centering
    \includegraphics[width=0.3\textwidth]{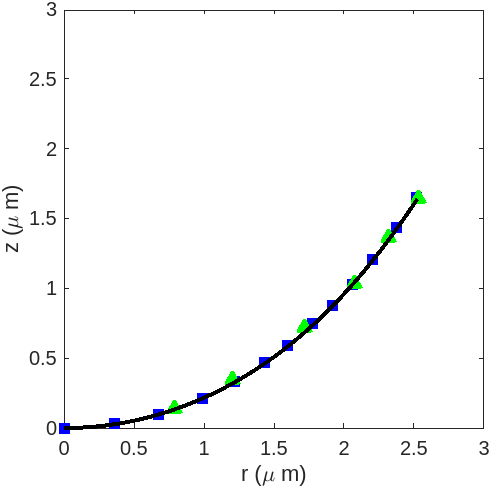} 
    \includegraphics[width=0.3\textwidth]{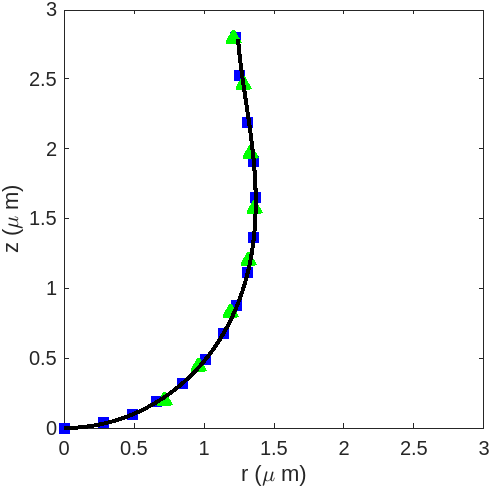} 
    \includegraphics[width=0.3\textwidth]{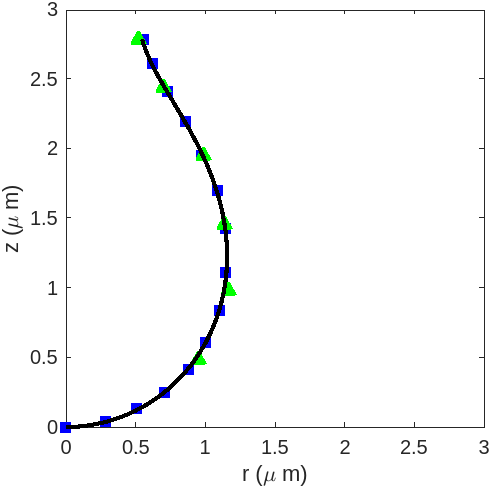} 
    \caption{Equilibrium shapes of an open membrane: the black solid curves are the numerical solutions obtained using our method; the blue squares depict the numerical results from \cite{Tu2010}; the green triangles show the experimental data from Fig.~3 (I to K) in \cite{Saitoh1998}. From left to right are cases $1,2$ and $3$ with reduced line tension  $\widetilde{\gamma} = 0.65, 0.78, 0.79\, (\mu m^{-1})$, and initial areas of $27.61, 23.15, 18.42\,(\mu m^{2})$ respectively.}
    \label{fig:equilibrium}
\end{figure}

We examine three cases for numerical comparison, where the spontaneous curvature and reduced Gaussian curvature rigidity are fixed as
$\gamma_g = -0.122$ and $c_0 = 0.2\, \mu m^{-1}$.
The negative value of $\gamma_g$ indicates that surfaces with positive Gaussian curvature are energetically more favorable than those with negative Gaussian curvature. The initial membrane shapes are selected to ensure their areas match the equilibrium areas reported in \cite{Tu2010}. The parameter values for each case are summarized in \Cref{tab:parameter_values}.
\begin{table}[htbp]
\centering
\caption{Parameter values for each case.}
\label{tab:parameter_values}
\begin{tabular}{|c|c|c|c|c|}
\hline
\textbf{Case} & $\gamma_g$ & $c_0\, (\mu m^{-1})$ & $\widetilde{\gamma}\, (\mu m^{-1})$ & $A_0\, (\mu m^2)$ \\
\hline
1 & -0.122 & 0.2 & 0.65 & 27.61 \\
2 & -0.122 & 0.2 & 0.78 & 23.15 \\
3 & -0.122 & 0.2 & 0.79 & 18.42 \\
\hline
\end{tabular}
\end{table}
The initial membrane shapes are selected to be close to the expected equilibrium configurations, as the membrane energy landscape can exhibit multiple equilibrium states. The dynamical equations are then solved over time until the total membrane energy ceases to decrease significantly. The computation is terminated when the relative change in total energy falls below $10^{-6}$. The time evolution of the total energy ${\cal E}$ for the three cases is shown in Figure~\ref{fig:e2-energy}.
The corresponding equilibrium shapes are presented in Figure~\ref{fig:equilibrium}. 

\subsection{Membrane dynamics from a spherical cap}
In this example, we investigate the dynamic behavior of an open membrane starting from an initial spherical cap shape. Previous numerical studies of membranes with a pore often assume that the membrane retains a spherical cap geometry throughout its evolution. This assumption simplifies the problem to a set of ODEs that govern a few geometric parameters \cite{Aubin2016,Ryham2018}. The justification for this approximation lies in the fact that a closed spherical membrane with constant mean curvature satisfies the equilibrium condition for normal stress balance. However, this assumption becomes less realistic in scenarios where the membrane has an open pore, particularly near the free edge \cite{Powers2002}.

Using our numerical method, we simulate the dynamics of an open membrane initialized as a spherical cap. The parameters are set as $\beta =  1$ and $\gamma_g = \gamma_l = H_0 = 0$, corresponding to the Willmore energy and zero line tension. The initial shape is defined by the generating curve:
\begin{equation}\label{eqn:eg3-init}
        X^r_0(s) = \sin(s), \quad
        X^z_0(s) = -\cos(s), \quad s \in [0, 0.9\pi].
\end{equation}

\begin{figure}[htbp]
    \centering
    \subfloat{\includegraphics[width=0.2\linewidth]{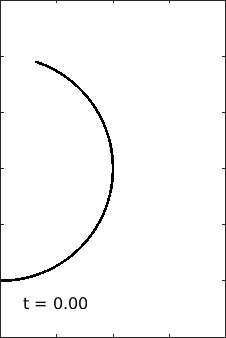}} \quad
    \subfloat{\includegraphics[width=0.2\linewidth]{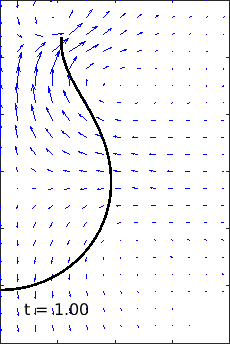}} \quad
    \subfloat{\includegraphics[width=0.2\linewidth]{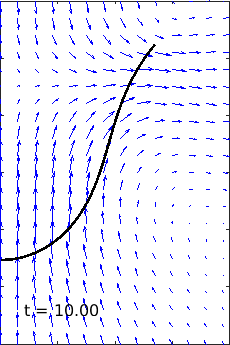}} \quad
    \subfloat{\includegraphics[width=0.2\linewidth]{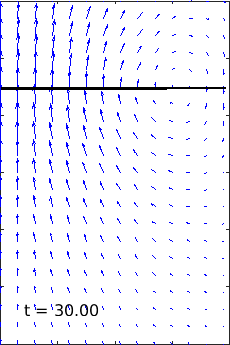}}\\[1ex]
    \subfloat{\includegraphics[width=0.2\linewidth]{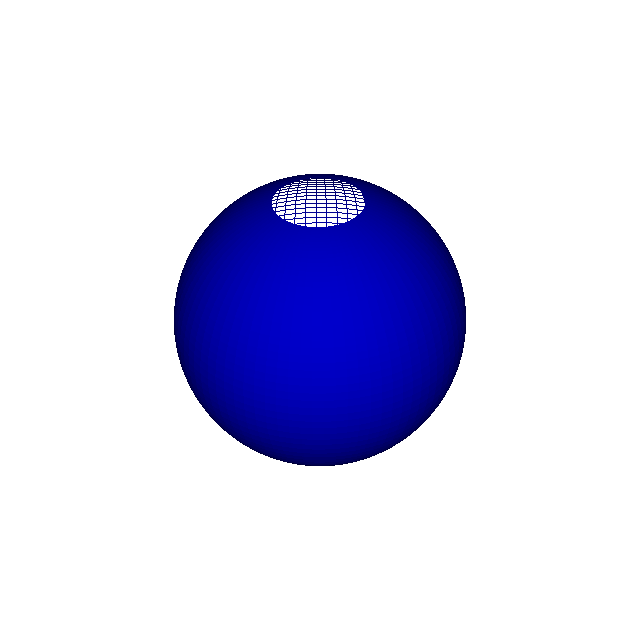}} \quad
    \subfloat{\includegraphics[width=0.2\linewidth]{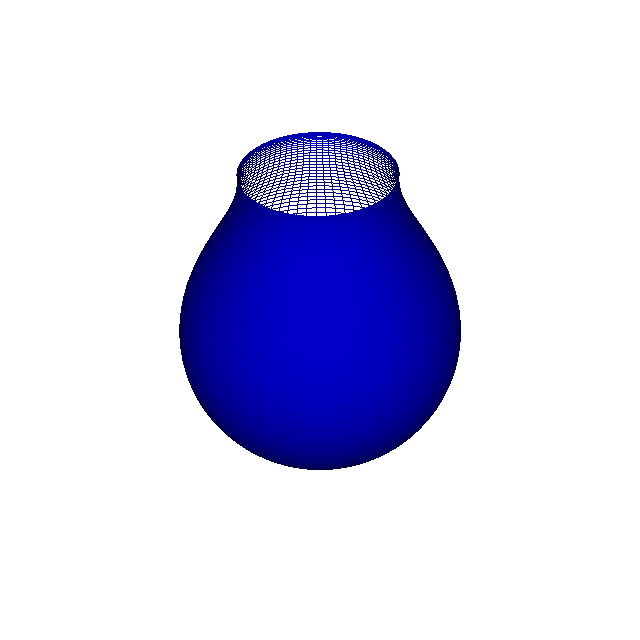}} \quad
    \subfloat{\includegraphics[width=0.2\linewidth]{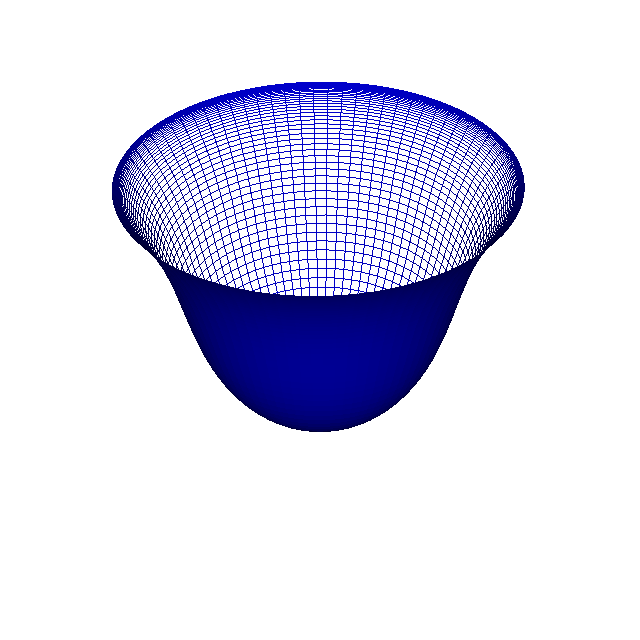}} \quad
    \subfloat{\includegraphics[width=0.2\linewidth]{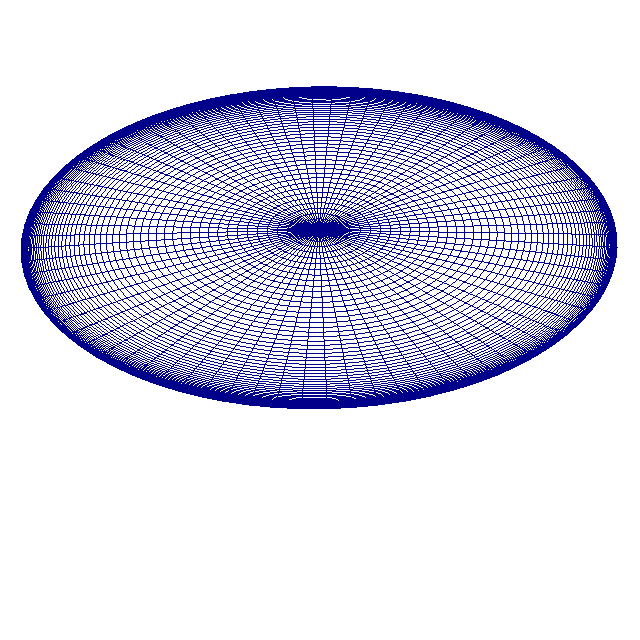}}
    \caption{Cross-sectional views (top row) and 3D visualizations (bottom row) of the membrane generating curve and the surrounding fluid velocity field at time points $t = 0$, $1$, $10$, and $30$.}
    \label{fig:eg3-a}
\end{figure}

\begin{figure}[htbp]
    \centering
    \includegraphics[width=0.35\linewidth]{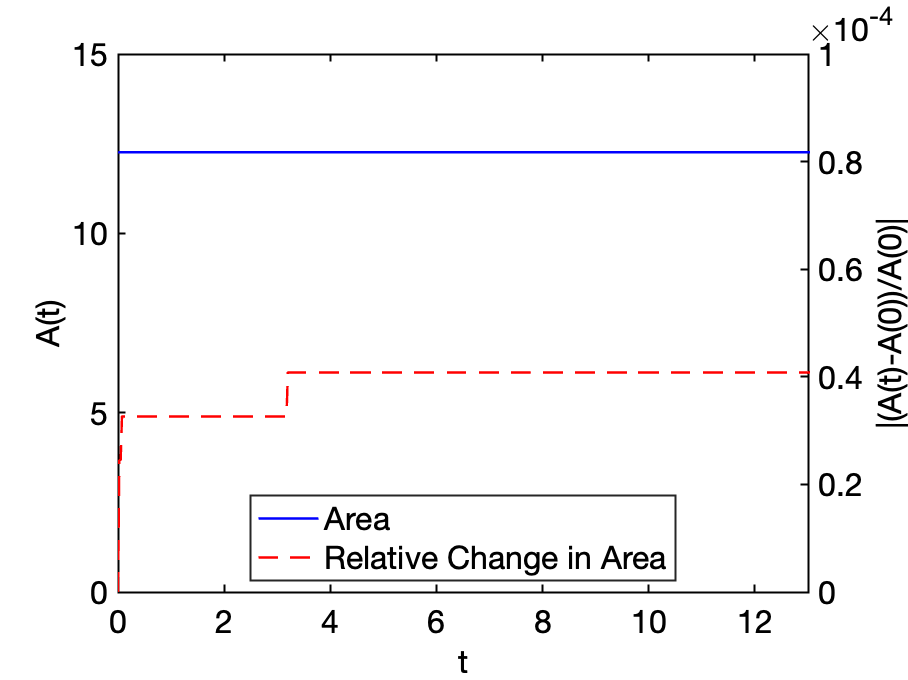}
    \includegraphics[width=0.35\linewidth]{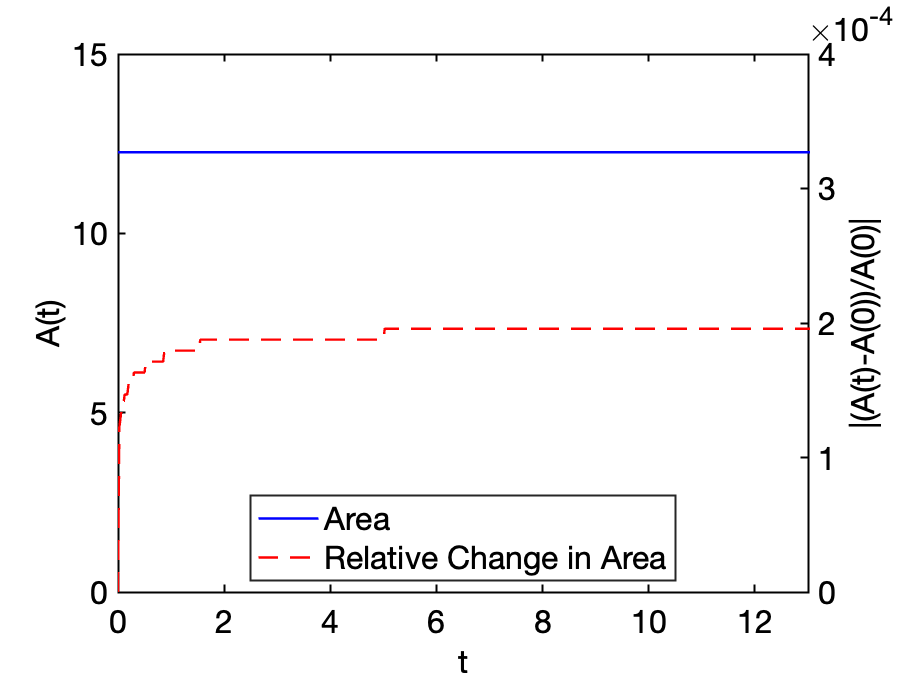}
    \caption{Time evolution of the membrane area and the corresponding relative change. Left: area profile for the flattening membrane; Right: area profile for the closing membrane.}
    \label{fig:eg3-area}
\end{figure}

The numerical results are shown in Figure~\ref{fig:eg3-a}. Initially, the membrane maintains a spherical cap shape. Over time, bending forces drive both bulk fluid flow and surface flow, causing the membrane to flatten. As the membrane evolves, deviations from the spherical cap geometry emerge, particularly near the open edge, where a neck-like structure begins to form.
This deviation can be understood through the boundary condition $H = 0$. While a spherical cap inherently possesses a non-zero constant mean curvature, this boundary condition is incompatible with the spherical cap assumption, leading to geometric deviations near the free edge.
We also plot the time evolution of the membrane area for the above two initial configurations in Figure~\ref{fig:eg3-area}. Since the membrane is locally inextensible, its total area should remain constant over time. Although the numerical solution introduces a slight variation due to discretization errors, the relative change in area remains below $0.05\%$ in both cases, demonstrating that the numerical method effectively preserves the area.

\begin{figure}[htbp]
    \centering
    \includegraphics[width=0.2\linewidth]{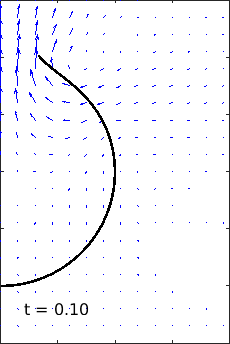}
    \includegraphics[width=0.2\linewidth]{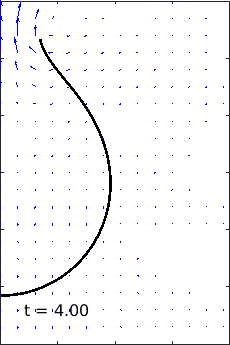}
    \includegraphics[width=0.2\linewidth]{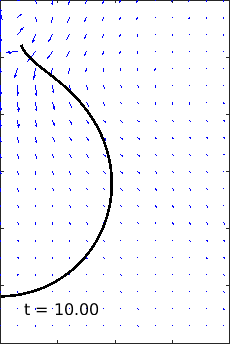}
    \includegraphics[width=0.2\linewidth]{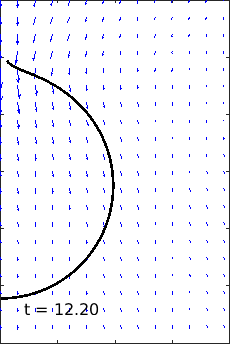}
    \caption{Cross-sectional views of the generating curve of the membrane and the velocity field of the fluid flow for $t=0.1,4,10,12.2$.}
    \label{fig:eg3-c}
\end{figure}

Next, we incorporate the influence of line tension by setting $\gamma_l = 0.5$. The line tension introduces an additional force along the open edge, acting to minimize its length and counteracting the bending force. The computation is terminated when the pore becomes sufficiently small to avoid numerical instabilities. As shown in Figure~\ref{fig:eg3-c}, a neck forms at an early stage due to the boundary condition on the mean curvature. Eventually, the pore closes while preserving the neck structure. In contrast, regions of the membrane far from the open edge maintain a nearly spherical cap shape.

\subsection{Boundary layer at the edge}
The above findings indicate that while the spherical cap assumption can provide an approximate description of the global membrane shape, particularly during early stages of pore formation, it fails to capture the detailed local behavior near the open edge. The imposed geometric boundary conditions play a crucial role in determining these local deviations.

\begin{figure}[htbp]
    \centering
    \includegraphics[width=0.3\linewidth]{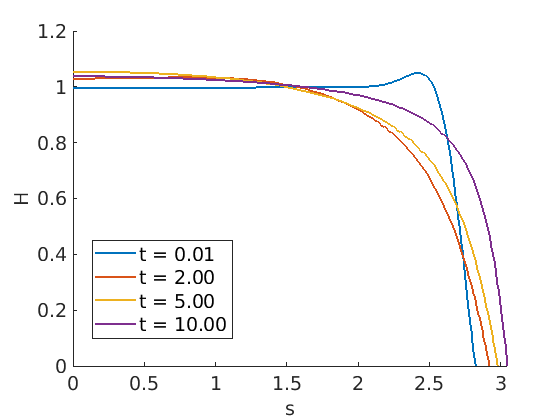}
    \includegraphics[width=0.3\linewidth]{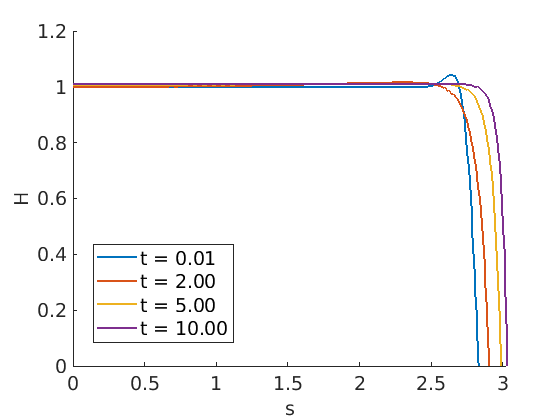}
    \includegraphics[width=0.3\linewidth]{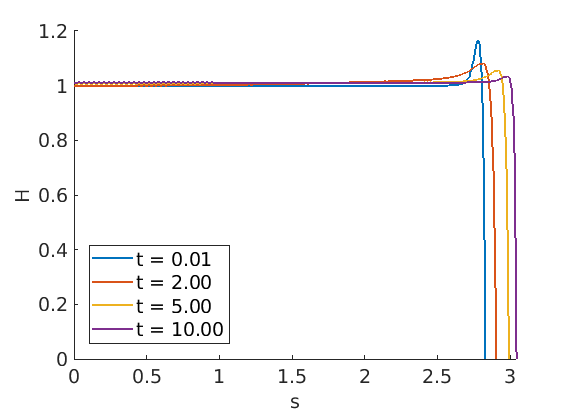}
    \caption{Numerical solution of the mean curvature on the membrane at different time  vs. arc-length parameter for $\gamma_l^{-1} = 1, 0.1, 0.01$ from left to right. The endpoint $s=0$ is the closed end in contacting with the $z-$axis, the other endpoint is the open edge.}
    \label{fig:eg3-H}
\end{figure}

As reported in \cite{Powers2002}, for small but non-zero $\epsilon=\gamma_l^{-1}$, an equilibrium membrane shape introduces an elastic boundary layer near the open edge. In our dynamic evolution equation, a similar boundary layer structure emerges during the membrane's evolution. To investigate this boundary layer more closely, we reduce $\gamma_l^{-1}$ from $1$ to $0.01$. Figure~\ref{fig:eg3-H} illustrates the mean curvature profiles for three different values of the ratio $\gamma_l^{-1} = 1, 0.1, 0.01$ and shows the variation of mean curvature along the arc-length parameter of the generating curve. The endpoint $s=0$ corresponds to the closed end in contact with the $z$-axis, while the other endpoint represents the open edge. 
The results reveal that as $\gamma_l^{-1}$ decreases, the boundary layer near the open edge becomes increasingly sharp and localized. For larger $\gamma_l^{-1}$, the mean curvature varies smoothly across the membrane, and the boundary layer effect is less pronounced. As $\gamma_l^{-1}$ becomes smaller, the curvature gradient steepens near the open edge, indicating a stronger influence of the elastic boundary layer.
This behavior suggests that the parameters $\gamma_l^{-1}$ play a critical role in controlling the membrane's local geometry near the open boundary. When $\gamma_l^{-1}$ is large, bending stiffness dominates, leading to a more uniform membrane profile. In contrast, when $\gamma_l^{-1}$ is small, the interplay between line tension and bending forces results in a highly localized boundary layer near the open edge.
These findings highlight the importance of accurately resolving the boundary layer structure in numerical simulations. Failure to capture this sharp transition can lead to significant inaccuracies in predicting the membrane's behavior near the open edge, particularly in scenarios where local geometric effects are dominant.

We compare our PDE model, which allows the membrane to adopt arbitrary shapes, with a simplified ODE model based on the spherical cap assumption. The ODE model assumes that the membrane maintains a spherical cap geometry throughout its evolution. This configuration is parameterized by three key variables: the sphere radius $R$, the cap angle $\alpha$, and the center position $z_0$. Using viscous force calculations from \cite{Ryham2018}, we can also derive an ODE system for the dynamics of $R$, $\alpha$, and $z_0$. 
Unlike the extensible membrane model in \cite{Ryham2018}, where surface tension is prescribed as a constant parameter, our ODE model enforces membrane inextensibility. In this framework, surface tension emerges as a Lagrange multiplier.

\begin{figure}[htbp]
    \centering
    \includegraphics[width=0.2\linewidth]{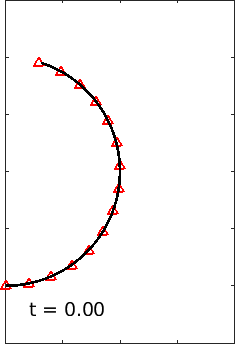}
    \includegraphics[width=0.2\linewidth]{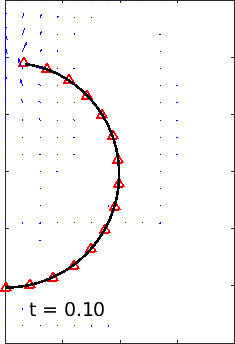}
    \includegraphics[width=0.2\linewidth]{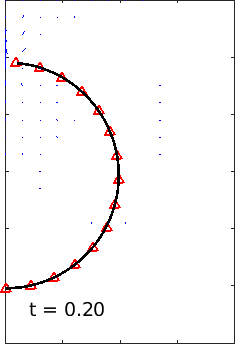}
    \includegraphics[width=0.2\linewidth]{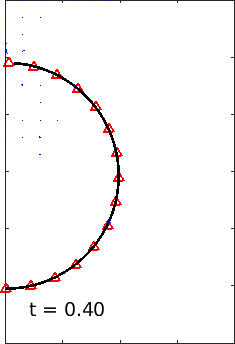}
    \caption{Numerical solutions for the open membrane obtained by solving the PDE model (solid dark line) and simplified ODE model (red triangles).}
    \label{fig:eg3-sph}
\end{figure}

For numerical simulations, we adopt physical parameters from \cite{Ryham2018}: an initial sphere radius of $20\,\mu$m, an initial hole radius of $6\,\mu$m, a line tension of $12\,\text{pN}$, fluid viscosity of $1\,\text{mP}\,\text{s}$, and membrane viscosity of $3 \times 10^{-9}\,\text{Pa}\,\text{m}\,\text{s}$. The characteristic spatial and time scales are chosen as $L = 20\,\mu$m and $T = 1$\,s.
Both the PDE and ODE models are initialized with the same shape described by \eqref{eqn:eg3-init}. The numerical solutions from both models are presented in Figure~\ref{fig:eg3-sph}. For such a parameter, the boundary layer at the edge is thin and {\color{black}the membrane predicted by the PDE} and ODE matches well in the outer region.

\begin{figure}[htbp]
    \centering
    \includegraphics[width=0.3\linewidth]{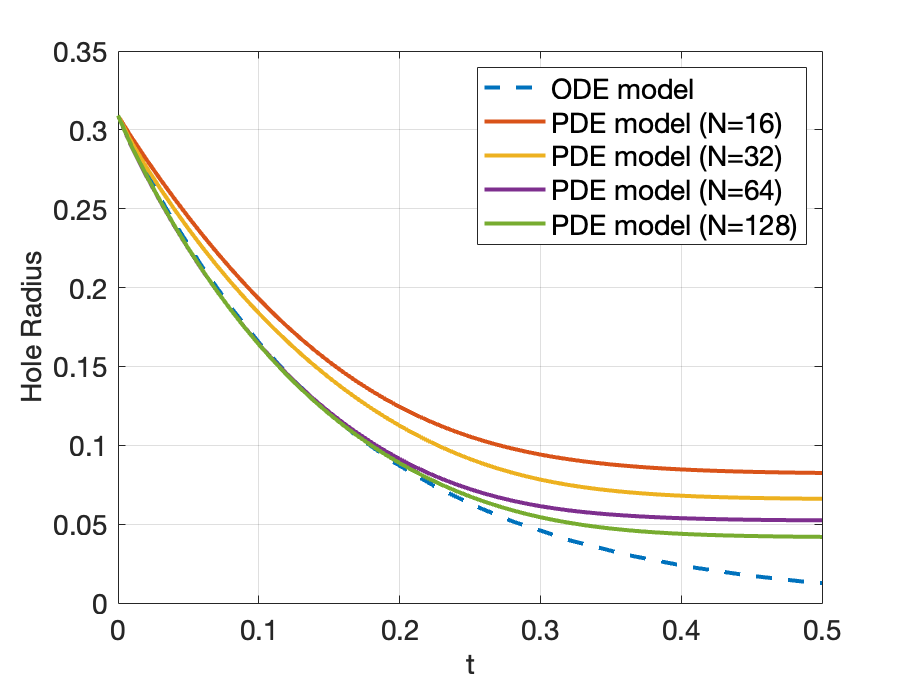}
    \includegraphics[width=0.3\linewidth]{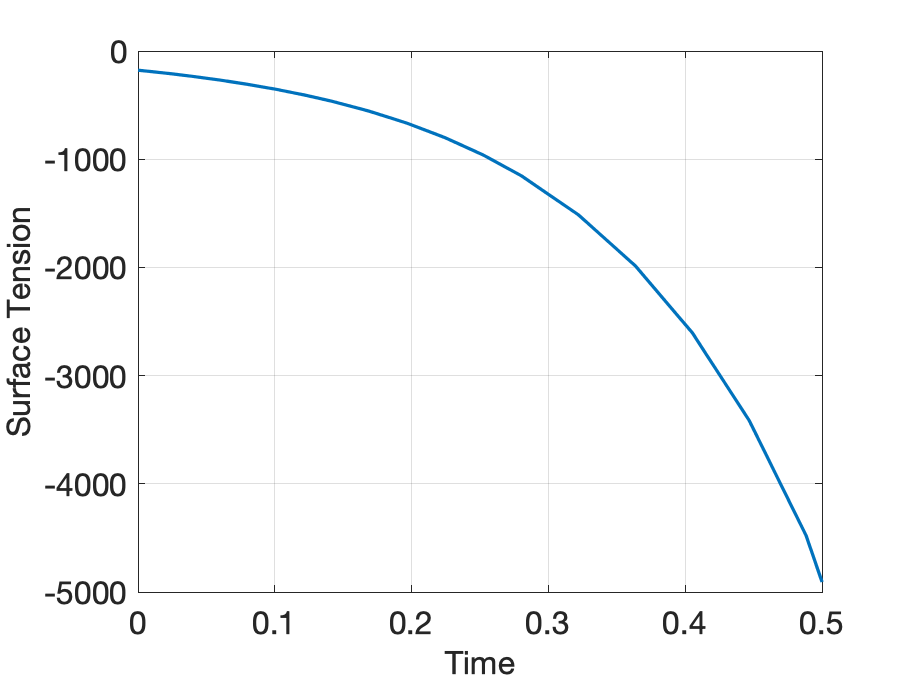}
    \includegraphics[width=0.3\linewidth]{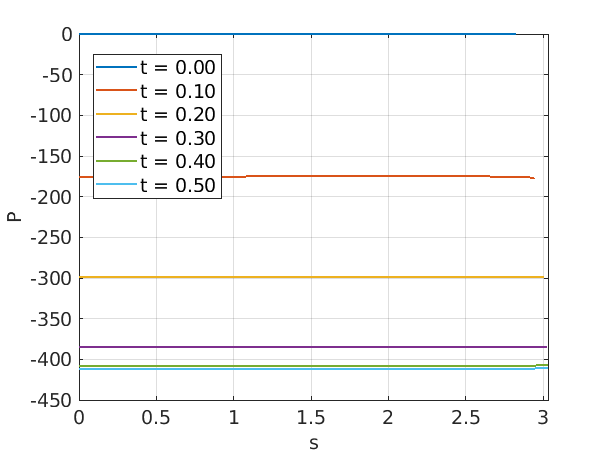}
    \caption{Left: hole radius of the spherical cap-shaped membrane predicted by the PDE and ODE models. Middle: surface tension predicted by the ODE model over time. Right: surface tension predicted by the PDE model, where the $x$-axis corresponds to the generating curve parameter.}
    \label{fig:eg3-sf}
\end{figure}

Figure~\ref{fig:eg3-sf} shows the evolution of the hole radius predicted by both the ODE and PDE models, with varying mesh refinement levels. At early stages ($t < 0.2$), the PDE model converges as the mesh is refined, and the ODE predictions closely match the PDE results. However, as the hole shrinks, discrepancies between the models become apparent. The ODE model predicts an exponential decay of the hole radius to zero, while the PDE model indicates that the hole radius approaches a steady state before further closure. 
In Figure~\ref{fig:eg3-sf}, we also compare the averaged surface tension predicted by the ODE model with the effective surface tension obtained from the PDE model. In the PDE simulation, the surface tension evolves toward a limiting value of approximately $-400$ and becomes nearly uniform across the membrane, indicating equilibrium. In contrast, the ODE model predicts a continuous decrease in surface tension, even when the hole becomes extremely small, with values significantly larger than those observed in the PDE model. This inconsistency highlights the inadequacy of the spherical cap assumption, particularly when considering the boundary layer effect of the membrane near the open edge.

\subsection{Membrane with different initial shapes}

In this example, we investigate the dynamics of lipid membranes starting from various initial shapes beyond the spherical cap configuration. We begin by examining a membrane initialized as a flat sheet to simulate the evolution of a closed vesicle originating from an open planar membrane. A flat sheet with zero spontaneous curvature represents a steady-state configuration under significant bending force and line tension. However, introducing a small degree of spontaneous curvature induces membrane bending, thereby amplifying the role of line tension in determining the membrane shape. For this study, the dimensionless parameters are set as $\gamma_g = 0$, $\beta = 1$, and $H_0 = 0.1$.

\begin{figure}[htbp]
    \centering
    \subfloat{\includegraphics[width=0.2\linewidth]{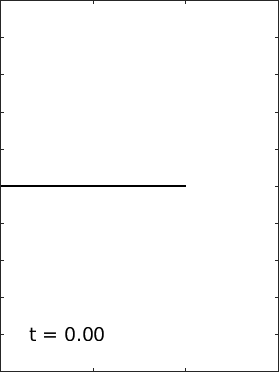}} \quad
    \subfloat{\includegraphics[width=0.2\linewidth]{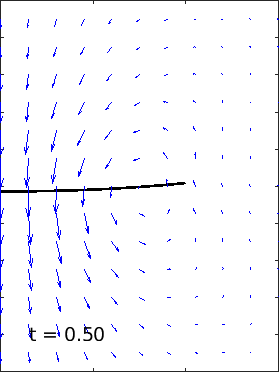}} \quad
    \subfloat{\includegraphics[width=0.2\linewidth]{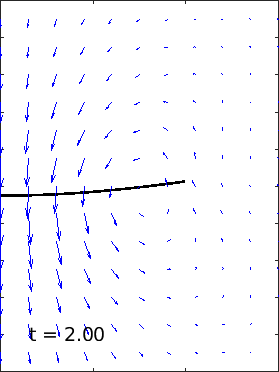}} \quad
    \subfloat{\includegraphics[width=0.2\linewidth]{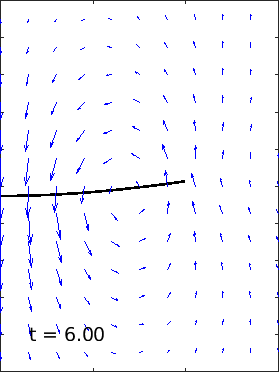}}\\[0.5ex]
    \subfloat{\includegraphics[width=0.2\linewidth]{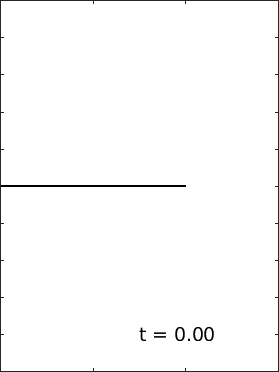}} \quad
    \subfloat{\includegraphics[width=0.2\linewidth]{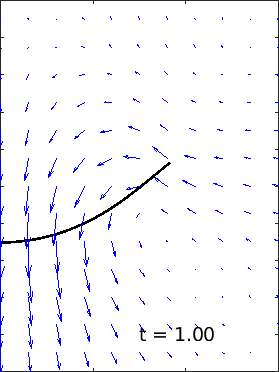}} \quad
    \subfloat{\includegraphics[width=0.2\linewidth]{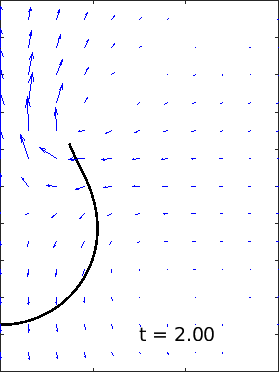}} \quad
    \subfloat{\includegraphics[width=0.2\linewidth]{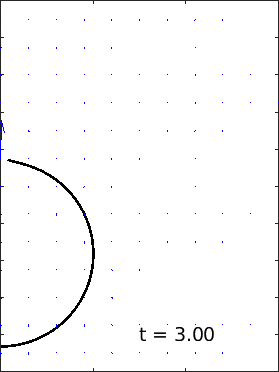}}
    \caption{Numerical solutions for an open membrane starting from a flat sheet. Top row: $\gamma_l = 1$ at time points $t = 0$, $0.5$, $2$, and $6$; bottom row: $\gamma_l = 5$ at time points $t = 0$, $1$, $2$, and $3$.}
    \label{fig:eg4-flat}
\end{figure}

First, we consider a small line tension value of $\gamma_l = 1$. Numerical results, shown in Figure~\ref{fig:eg4-flat}, indicate that the small spontaneous curvature causes the membrane to bend slightly, approaching a near-equilibrium configuration. However, the low line tension is insufficient to drive the membrane toward full closure into a vesicle.

Subsequently, we increase the line tension to $\gamma_l = 5$. As illustrated in Figure~\ref{fig:eg4-flat}, the elevated line tension actively drives the open edge of the membrane to close, forming a nearly closed vesicle structure.
To further investigate membrane dynamics, we consider an initial membrane shape resembling a biconcave red blood cell with a small pore. The parametric form of this shape is given by
\begin{equation}
\begin{aligned}
    &r(\eta) = 2\sqrt{\eta(1-\eta)},\quad z(\eta) = 0.7(2\eta-1)-0.6(2\eta-1)^3 +0.05(2\eta-1)^5, \\
    &\eta(s) = 0.475(1-\cos(\pi s)),\quad s\in [0,1].
\end{aligned}
\end{equation}

\begin{figure}[htbp]
    \centering
    \subfloat{\includegraphics[width=0.2\linewidth]{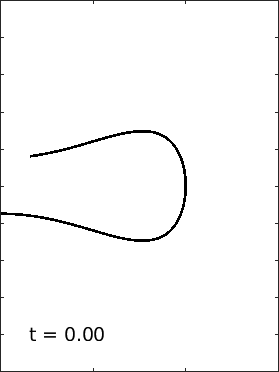}} \quad
    \subfloat{\includegraphics[width=0.2\linewidth]{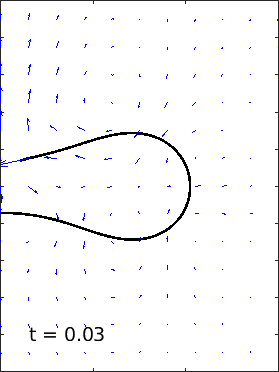}} \quad
    \subfloat{\includegraphics[width=0.2\linewidth]{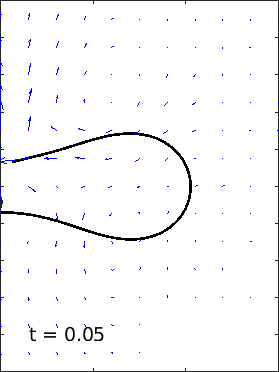}} \quad
    \subfloat{\includegraphics[width=0.2\linewidth]{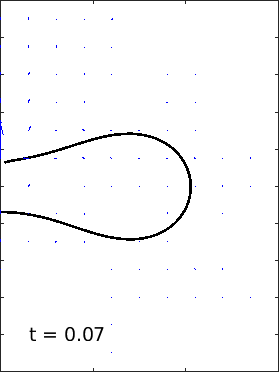}}\\[0.5ex]
    \subfloat{\includegraphics[width=0.2\linewidth]{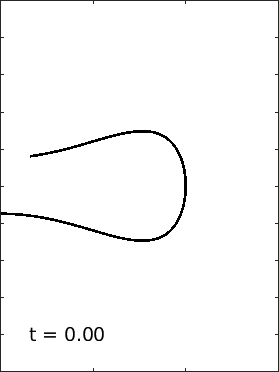}} \quad
    \subfloat{\includegraphics[width=0.2\linewidth]{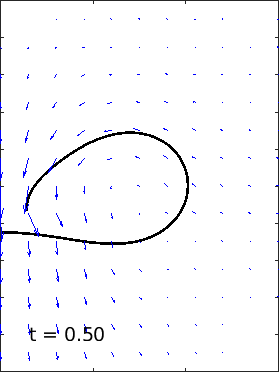}} \quad
    \subfloat{\includegraphics[width=0.2\linewidth]{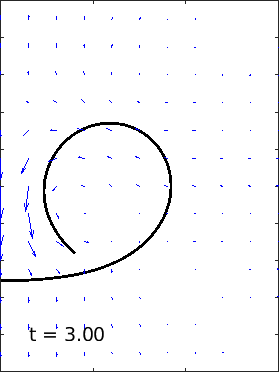}} \quad
    \subfloat{\includegraphics[width=0.2\linewidth]{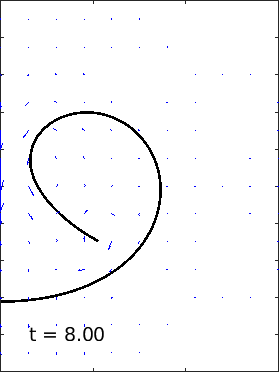}}
    \caption{Numerical solutions for an open membrane initially shaped as a biconcave red blood cell with a small pore, for $\gamma_l = 5$ (top row) and $\gamma_l = 1$ (bottom row).}
    \label{fig:eg4-bic}
\end{figure}

Using the same dimensionless parameters as before ($\gamma_g = 0$, $\beta = 1$, and $H_0 = 0.1$), we first consider a strong line tension ($\gamma_l = 5$). As shown in Figure~\ref{fig:eg4-bic}, the strong line tension drives the small pore to close rapidly, leaving the enclosed volume of the cell essentially unchanged, meaning {\color{black}the resulting steady state is not affected much}. This scenario suggests efficient healing of membrane damage without significant deformation.

Conversely, when the line tension is reduced to $\gamma_l = 1$, the pore does not close immediately. Instead, the membrane begins to deform into an alternative configuration, preventing it from returning to a closed vesicle state. This behavior, illustrated in Figure~\ref{fig:eg4-bic}, underscores the critical role of line tension in the formation, stability, and repair processes of lipid membranes.

\section{Disscusion}\label{sec:6}

For an open, inextensible membrane immersed in Stokes flow, we derive a continuum model based on the principle of maximum dissipation within an energy variational framework. The membrane is represented as a zero-thickness open surface endowed with a membrane energy comprising the Helfrich bending energy and a line tension along its open edge. Dissipative effects from both the bulk and membrane viscosities are incorporated via Rayleigh dissipation functionals, while incompressibility of the surrounding fluid and inextensibility of the membrane are enforced through Lagrange multipliers. Taking first variations and balancing surface and boundary forces yields the coupled PDEs governing membrane–fluid interactions.

To conduct numerical simulations for the PDE model, we consider the axisymmetric case for both the membrane and fluid flow. For this scenario, we adopt a boundary integral formulation, reducing the 3D problem into a 1D problem along the generating curve. Due to the complex coupling between the boundary integral and geometric quantities, we developed a hybrid BEM-FEM method by rewriting the high-order PDE into a weak mixed form, which only contains first derivatives. The time discretization is semi-implicit: low-order derivatives, nonlinear terms and boundary terms are treated explicitly, while high-order derivatives, linear terms are treated implicitly. This ensures that the final scheme is both stable and can be solved efficiently. The single-layer density typically exhibits singular behavior near the open edge. To address this, we designed a local mesh refinement strategy that increases the resolution of the computational mesh near the open edge, improving the scheme's ability to capture singularities. Using the developed axisymmetric numerical scheme, we applied it to model open membranes in different cases. The present numerical method can be applied to investigate the behavior of open membranes interacting with the bulk flow, such as the relaxation dynamics of membranes under the influence of the bulk fluid and the dynamics of pores in membranes influenced by fluid leakage.

Lipid bilayer membranes are nearly inextensible, and a weak area compressibility may result in fluid leakage that alters both the membrane dynamics \cite{fu2022two} and the pore dynamics \cite{Ryham2018}. How does such weak area compressibility affect the transition from a (closed) vesicle to an open membrane? To answer this question we have been working on improving our PDE model for the lipid bilayer membrane to account for the physical area compressibility within the boundary integral framework. 
{\color{black}Beyond lipid membranes, our framework may also apply to colloidal membranes and help interpret experimentally observed colloidosome assembly and disassembly pathways \cite{Adkins2025}.}

Furthermore the current numerical method can be improved in several directions. By carefully studying the singularity behavior of the solution near the open edge, it should be possible to develop a high-order numerical method that incorporates an analytical expression for the singularity structure, reducing the need for mesh refinement. Developing numerical schemes to perform full 3D simulations for membranes with general shapes is another challenging but interesting direction.


\bibliographystyle{siamplain}
\bibliography{references}
\end{document}

%% file: ex_shared.tex
\usepackage{lipsum}
\usepackage{amsfonts}
\usepackage{graphicx}
\usepackage{epstopdf}
\usepackage{algorithmic}
\ifpdf
  \DeclareGraphicsExtensions{.eps,.pdf,.png,.jpg}
\else
  \DeclareGraphicsExtensions{.eps}
\fi

\newsiamremark{remark}{Remark}
\newsiamremark{hypothesis}{Hypothesis}
\crefname{hypothesis}{Hypothesis}{Hypotheses}
\newsiamthm{claim}{Claim}

\headers{Open Membrane in Stokes Flow}{H. Zhou, Y, Young, and Y. Mori}

\title{Modeling and Simulation of Open Membranes in Stokes Flow with Mixed-Dimensional Coupling\thanks{Submitted to the editors DATE.
\funding{YM was partially supported by the National Science Foundation (NSF) grant DMR2309034 (USA), and the Math+X award from the Simons Foundation (Award ID 234606). YNY acknowledges support from the NSF under Award No.
DMS-1951600 and DMS-251071, and from the Flatiron Institute, part of Simons Foundation.}}}

\author{Han Zhou\thanks{Department of Mathematics, University of Pennsylvania, Philadelphia, PA 19104, USA 
  (\email{hzhou24@sas.upenn.edu}).}
\and Yuan-Nan Young\thanks{Department of Mathematical Sciences, New Jersey Institute of Technology, Newark, NJ 07102, USA 
  (\email{yyoung@njit.edu}).}
\and Yoichiro Mori\thanks{Department of Mathematics, University of Pennsylvania, Philadelphia, PA 19104, USA 
  (\email{y1mori@sas.upenn.edu}).}
}

\usepackage{amsopn}

\makeatletter
\newcommand*{\addFileDependency}[1]{
  \typeout{(#1)}
  \@addtofilelist{#1}
  \IfFileExists{#1}{}{\typeout{No file #1.}}
}
\makeatother

\newcommand*{\myexternaldocument}[1]{%
    \externaldocument{#1}%
    \addFileDependency{#1.tex}%
    \addFileDependency{#1.aux}%
}

%% file: references.bib
@article{hyon2010energetic,
  title={Energetic variational approach in complex fluids: maximum dissipation principle},
  author={Hyon, Yunkyong and Kwak, Do Young and Liu, Chun},
  journal={Discrete Contin. Dyn. Syst},
  volume={26},
  number={4},
  pages={1291--1304},
  year={2010}
}

@article{seifert1997configurations,
  title={Configurations of fluid membranes and vesicles},
  author={Seifert, Udo},
  journal={Advances in physics},
  volume={46},
  number={1},
  pages={13--137},
  year={1997},
  publisher={Taylor \& Francis}
}

@article{hu2015hybrid,
  title={A hybrid immersed boundary and immersed interface method for electrohydrodynamic simulations},
  author={Hu, Wei-Fan and Lai, Ming-Chih and Young, Yuan-Nan},
  journal={Journal of Computational Physics},
  volume={282},
  pages={47--61},
  year={2015},
  publisher={Elsevier}
}

@article{hu2016vesicle,
  title={Vesicle electrohydrodynamic simulations by coupling immersed boundary and immersed interface method},
  author={Hu, Wei-Fan and Lai, Ming-Chih and Seol, Yunchang and Young, Yuan-Nan},
  journal={Journal of Computational Physics},
  volume={317},
  pages={66--81},
  year={2016},
  publisher={Elsevier}
}

@article{fu2020simulation,
  title={Simulation of multiscale hydrophobic lipid dynamics via efficient integral equation methods},
  author={Fu, Szu-Pei P and Ryham, Rolf and Kl\"{o}ckner, Andreas and Wala, Matt and Jiang, Shidong and Young, Yuan-Nan},
  journal={Multiscale Modeling \& Simulation},
  volume={18},
  number={1},
  pages={79--103},
  year={2020},
  publisher={SIAM}
}

@article{ryham2011aqueous,
  title={Aqueous viscosity is the primary source of friction in lipidic pore dynamics},
  author={Ryham, Rolf and Berezovik, Irina and Cohen, Fredric S},
  journal={Biophysical Journal},
  volume={101},
  number={12},
  pages={2929--2938},
  year={2011},
  publisher={Elsevier}
}

@article{ryham2013teardrop,
  title={Teardrop shapes minimize bending energy of fusion pores connecting planar bilayers},
  author={Ryham, Rolf J and Ward, Mark A and Cohen, Fredric S},
  journal={Physical Review E—Statistical, Nonlinear, and Soft Matter Physics},
  volume={88},
  number={6},
  pages={062701},
  year={2013},
  publisher={APS}
}

@article{Powers2002,
   author = {Thomas R. Powers and Greg Huber and Raymond E. Goldstein},
   issn = {1063651X},
   issue = {4},
   journal = {Physical Review E - Statistical Physics, Plasmas, Fluids, and Related Interdisciplinary Topics},
   pages = {11},
   pmid = {12005867},
   title = {Fluid-membrane tethers: Minimal surfaces and elastic boundary layers},
   volume = {65},
   year = {2002},
}

@article{Chuan2021,
author = {Chuan, Kian and Lai, Ming-chih and Seol, Yunchang},
issn = {0021-9991},
journal = {Journal of Computational Physics},
pages = {110090},
publisher = {Elsevier Inc.},
title = {{An immersed boundary projection method for incompressible interface simulations in 3D flows}},
volume = {430},
year = {2021}
}

@article{Reuther2016,
author = {Reuther, Sebastian and Voigt, Axel},
issn = {0021-9991},
journal = {Journal of Computational Physics},
pages = {850--858},
publisher = {Elsevier Inc.},
title = {{Incompressible two-phase flows with an inextensible Newtonian fluid interface}},
volume = {322},
year = {2016}
}

@article{Kim2010,
author = {Kim, Yongsam and Lai, Ming-chih},
doi = {10.1016/j.jcp.2010.03.020},
issn = {0021-9991},
journal = {Journal of Computational Physics},
number = {12},
pages = {4840--4853},
publisher = {Elsevier Inc.},
title = {{Simulating the dynamics of inextensible vesicles by the penalty immersed boundary method}},
url = {http://dx.doi.org/10.1016/j.jcp.2010.03.020},
volume = {229},
year = {2010}
}

@article{Farutin2014,
author = {Farutin, Alexander and Biben, Thierry and Misbah, Chaouqi},
doi = {10.1016/j.jcp.2014.07.008},
issn = {0021-9991},
journal = {Journal of Computational Physics},
pages = {539--568},
publisher = {Elsevier Inc.},
title = {{3D numerical simulations of vesicle and inextensible capsule dynamics}},
url = {http://dx.doi.org/10.1016/j.jcp.2014.07.008},
volume = {275},
year = {2014}
}

@article{Veerapaneni2009,
author = {Veerapaneni, Shravan K. and Gueyffier, Denis and Biros, George and Zorin, Denis},
doi = {10.1016/j.jcp.2009.06.020},
isbn = {2129983510},
issn = {00219991},
journal = {Journal of Computational Physics},
month = {oct},
number = {19},
pages = {7233--7249},
publisher = {Elsevier Inc.},
title = {{A numerical method for simulating the dynamics of 3D axisymmetric vesicles suspended in viscous flows}},
url = {http://dx.doi.org/10.1016/j.jcp.2009.06.020 https://linkinghub.elsevier.com/retrieve/pii/S0021999109003441},
volume = {228},
year = {2009}
}

@article{Lai2019,
author = {Lai, Ming Chih and Ong, Kian Chuan},
doi = {10.1137/18M1210277},
issn = {10957197},
journal = {SIAM Journal on Scientific Computing},
month = {jan},
number = {4},
pages = {B649--B668},
title = {{Unconditionally energy stable schemes for the inextensible interface problem with bending}},
url = {https://epubs.siam.org/doi/10.1137/18M1210277},
volume = {41},
year = {2019}
}

@incollection{Barrett2020,
author = {Barrett, John W. and Garcke, Harald and N{\"{u}}rnberg, Robert},
booktitle = {Handbook of Numerical Analysis},
doi = {10.1016/bs.hna.2019.05.002},
eprint = {1903.09462},
isbn = {9780444640031},
issn = {15708659},
pages = {275--423},
title = {{Parametric finite element approximations of curvature-driven interface evolutions}},
url = {https://linkinghub.elsevier.com/retrieve/pii/S1570865919300055},
volume = {21},
year = {2020}
}

@article{Jia2022,
archivePrefix = {arXiv},
arxivId = {2201.00908},
author = {Jia, Leroy L. and Shelley, Michael J.},
doi = {10.1017/jfm.2022.550},
eprint = {2201.00908},
issn = {0022-1120},
journal = {Journal of Fluid Mechanics},
month = {oct},
pages = {A1},
title = {{The role of monolayer viscosity in Langmuir film hole closure dynamics}},
url = {https://www.cambridge.org/core/product/identifier/S002211202200550X/type/journal_article},
volume = {948},
year = {2022}
}

@article{Ryham2018,
author = {Ryham, Rolf J.},
doi = {10.1017/jfm.2017.807},
issn = {0022-1120},
journal = {Journal of Fluid Mechanics},
month = {feb},
pages = {502--531},
title = {{On the viscous flows of leak-out and spherical cap natation}},
url = {https://www.cambridge.org/core/product/identifier/S0022112017008072/type/journal_article},
volume = {836},
year = {2018}
}

@article{Alpert1999,
author = {Alpert, Bradley K},
doi = {10.1137/S1064827597325141},
issn = {1064-8275},
journal = {SIAM Journal on Scientific Computing},
month = {jan},
number = {5},
pages = {1551--1584},
title = {{Hybrid Gauss-Trapezoidal Quadrature Rules}},
url = {http://epubs.siam.org/doi/10.1137/S1064827597325141},
volume = {20},
year = {1999}
}

@manual{eigenweb,
  title  = {Eigen: A C++ template library for linear algebra},
  author = {Ga\"{e}l Guennebaud and Beno\^{i}t Jacob and others},
  year   = {2010},
  url    = {http://eigen.tuxfamily.org}
}

@article{Saitoh1998,
author = {Saitoh, Akihiko and Takiguchi, Kingo and Tanaka, Yohko and Hotani, Hirokazu},
doi = {10.1073/pnas.95.3.1026},
issn = {0027-8424},
journal = {Proceedings of the National Academy of Sciences},
month = {feb},
number = {3},
pages = {1026--1031},
pmid = {9448279},
title = {{Opening-up of liposomal membranes by talin}},
url = {https://pnas.org/doi/full/10.1073/pnas.95.3.1026},
volume = {95},
year = {1998}
}

@article{balchunas2019equation,
  title={Equation of state of colloidal membranes},
  author={Balchunas, Andrew J and Cabanas, Rafael A and Zakhary, Mark J and Gibaud, Thomas and Fraden, Seth and Sharma, Prerna and Hagan, Michael F and Dogic, Zvonimir},
  journal={Soft Matter},
  volume={15},
  number={34},
  pages={6791--6802},
  year={2019},
  publisher={Royal Society of Chemistry}
}

@article{balchunas2020force,
  title={Force-induced formation of twisted chiral ribbons},
  author={Balchunas, Andrew and Jia, Leroy L and Zakhary, Mark J and Robaszewski, Joanna and Gibaud, Thomas and Dogic, Zvonimir and Pelcovits, Robert A and Powers, Thomas R},
  journal={Physical Review Letters},
  volume={125},
  number={1},
  pages={018002},
  year={2020},
  publisher={APS}
}

@article{Capovilla2002,
   author = {R. Capovilla and J. Guven and J. A. Santiago},
   doi = {10.1103/PhysRevE.66.021607},
   issn = {1063651X},
   issue = {2},
   journal = {Physical Review E - Statistical Physics, Plasmas, Fluids, and Related Interdisciplinary Topics},
   pages = {1-7},
   title = {Lipid membranes with an edge},
   volume = {66},
   year = {2002},
}

@article{Klenow2024,
   author = {Martin Berg Klenow and Magnus Staal Vigsø and Weria Pezeshkian and Jesper Nylandsted and Michael Andersen Lomholt and Adam Cohen Simonsen},
   doi = {10.1016/j.bpj.2024.05.027},
   issn = {00063495},
   issue = {13},
   journal = {Biophysical Journal},
   month = {7},
   pages = {1827-1837},
   publisher = {Biophysical Society},
   title = {Shape of the membrane neck around a hole during plasma membrane repair},
   volume = {123},
   url = {https://doi.org/10.1016/j.bpj.2024.05.027 https://linkinghub.elsevier.com/retrieve/pii/S0006349524003746},
   year = {2024},
}

@article{Palmer2024,
   author = {Bennett Palmer and Álvaro Pámpano},
   doi = {10.1016/j.na.2023.113393},
   issn = {0362546X},
   journal = {Nonlinear Analysis, Theory, Methods and Applications},
   pages = {113393},
   publisher = {Elsevier Ltd},
   title = {Symmetry breaking bifurcation of membranes with boundary},
   volume = {238},
   url = {https://doi.org/10.1016/j.na.2023.113393},
   year = {2024},
}

@article{Malik2022,
   author = {Vinit Kumar Malik and On Shun Pak and Jie Feng},
   doi = {10.1103/PhysRevApplied.17.024032},
   issn = {23317019},
   issue = {2},
   journal = {Physical Review Applied},
   pages = {1},
   publisher = {American Physical Society},
   title = {Pore Dynamics of Lipid Vesicles under Light-Induced Osmotic Stress},
   volume = {17},
   url = {https://doi.org/10.1103/PhysRevApplied.17.024032},
   year = {2022},
}

@article{Palmer2022,
   author = {Bennett Palmer and Álvaro Pámpano},
   doi = {10.1007/s00526-022-02188-6},
   issn = {14320835},
   issue = {3},
   journal = {Calculus of Variations and Partial Differential Equations},
   pages = {1-28},
   publisher = {Springer Berlin Heidelberg},
   title = {The Euler–Helfrich functional},
   volume = {61},
   url = {https://doi.org/10.1007/s00526-022-02188-6},
   year = {2022},
}

@article{Barrett2021,
   author = {John W. Barrett and Harald Garcke and Robert Nürnberg},
   doi = {10.1051/m2an/2021014},
   issn = {0764-583X},
   issue = {3},
   journal = {ESAIM: Mathematical Modelling and Numerical Analysis},
   month = {5},
   pages = {833-885},
   title = {Stable approximations for axisymmetric Willmore flow for closed and open surfaces},
   volume = {55},
   url = {https://www.esaim-m2an.org/10.1051/m2an/2021014},
   year = {2021},
}

@article{Palmer2021,
   author = {Bennett Palmer and Álvaro Pámpano},
   doi = {10.1007/s00332-021-09679-4},
   isbn = {0033202109},
   issn = {0938-8974},
   issue = {1},
   journal = {Journal of Nonlinear Science},
   month = {2},
   pages = {23},
   publisher = {Springer US},
   title = {Minimizing Configurations for Elastic Surface Energies with Elastic Boundaries},
   volume = {31},
   url = {https://doi.org/10.1007/s00332-021-09679-4 http://link.springer.com/10.1007/s00332-021-09679-4},
   year = {2021},
}

@article{Zhou2018,
   author = {Xiaohua Zhou},
   doi = {10.1016/j.ijnonlinmec.2018.08.019},
   issn = {00207462},
   issue = {May},
   journal = {International Journal of Non-Linear Mechanics},
   pages = {25-28},
   publisher = {Elsevier Ltd},
   title = {An integral case of the axisymmetric shape equation of open vesicles with free edges},
   volume = {106},
   url = {https://doi.org/10.1016/j.ijnonlinmec.2018.08.019},
   year = {2018},
}

@article{Barrett2017,
   author = {John W. Barrett and Harald Garcke and Robert Nürnberg},
   doi = {10.1093/imanum/drx006},
   issn = {0272-4979},
   issue = {4},
   journal = {IMA Journal of Numerical Analysis},
   month = {3},
   pages = {1657-1709},
   title = {Stable variational approximations of boundary value problems for Willmore flow with Gaussian curvature},
   volume = {37},
   url = {https://academic.oup.com/imajna/article-lookup/doi/10.1093/imanum/drx006},
   year = {2017},
}

@article{Tang2017,
   author = {Sindy K Y Tang and Wallace F Marshall},
   doi = {10.1126/science.aam6496},
   issn = {0036-8075},
   issue = {6342},
   journal = {Science},
   month = {6},
   pages = {1022-1025},
   title = {Self-repairing cells: How single cells heal membrane ruptures and restore lost structures},
   volume = {356},
   url = {https://www.science.org/doi/10.1126/science.aam6496},
   year = {2017},
}

@article{Aubin2016,
   author = {Christopher A. Aubin and Rolf J. Ryham},
   doi = {10.1017/jfm.2015.699},
   issn = {0022-1120},
   journal = {Journal of Fluid Mechanics},
   month = {2},
   pages = {228-245},
   title = {Stokes flow for a shrinking pore},
   volume = {788},
   url = {https://www.cambridge.org/core/product/identifier/S0022112015006990/type/journal_article},
   year = {2016},
}

@article{Deserno2015,
   author = {Markus Deserno},
   doi = {10.1016/j.chemphyslip.2014.05.001},
   issn = {00093084},
   journal = {Chemistry and Physics of Lipids},
   month = {1},
   pages = {11-45},
   pmid = {24835737},
   publisher = {Elsevier Ireland Ltd},
   title = {Fluid lipid membranes: From differential geometry to curvature stresses},
   volume = {185},
   url = {http://dx.doi.org/10.1016/j.chemphyslip.2014.05.001 https://linkinghub.elsevier.com/retrieve/pii/S000930841400053X},
   year = {2015},
}

@article{Cooper2015,
   author = {Sandra T. Cooper and Paul L. McNeil},
   doi = {10.1152/physrev.00037.2014},
   issn = {0031-9333},
   issue = {4},
   journal = {Physiological Reviews},
   month = {10},
   pages = {1205-1240},
   title = {Membrane Repair: Mechanisms and Pathophysiology},
   volume = {95},
   url = {https://www.physiology.org/doi/10.1152/physrev.00037.2014},
   year = {2015},
}

@article{Asgari2015,
   author = {Meisam Asgari and Aisa Biria},
   doi = {10.1016/j.ijnonlinmec.2015.06.001},
   issn = {00207462},
   journal = {International Journal of Non-Linear Mechanics},
   pages = {135-143},
   publisher = {Elsevier},
   title = {Free energy of the edge of an open lipid bilayer based on the interactions of its constituent molecules},
   volume = {76},
   url = {http://dx.doi.org/10.1016/j.ijnonlinmec.2015.06.001},
   year = {2015},
}

@article{fu2022two,
  title={Two-dimensional hydrodynamics of a Janus particle vesicle},
  author={Fu, Szu-Pei and Quaife, Bryan and Ryham, Rolf and Young, Y-N},
  journal={Journal of Fluid Mechanics},
  volume={941},
  pages={A41},
  year={2022},
  publisher={Cambridge University Press}
}

@inbook{Biria2013,
   author = {Aisa Biria and Mohsen Maleki and Eliot Fried},
   doi = {10.1016/B978-0-12-396522-6.00001-3},
   edition = {1},
   isbn = {9780123965226},
   issn = {00652156},
   journal = {Advances in Applied Mechanics},
   pages = {1-68},
   publisher = {Elsevier Inc.},
   title = {Continuum Theory for the Edge of an Open Lipid Bilayer},
   volume = {46},
   url = {http://dx.doi.org/10.1016/B978-0-12-396522-6.00001-3 https://linkinghub.elsevier.com/retrieve/pii/B9780123965226000013},
   year = {2013},
}

@article{Tu2011,
   author = {Zhanchun Tu},
   doi = {10.7546/jgsp-24-2011-45-75},
   issn = {13145673},
   journal = {Journal of Geometry and Symmetry in Physics},
   pages = {45-75},
   title = {Geometry of membranes},
   volume = {24},
   year = {2011},
}

@article{Laadhari2010,
   author = {A Laadhari and C Misbah and P Saramito},
   doi = {10.1016/j.physd.2010.04.001},
   issn = {0167-2789},
   issue = {16},
   journal = {Physica D},
   pages = {1567-1572},
   publisher = {Elsevier B.V.},
   title = {On the equilibrium equation for a generalized biological membrane energy by using a shape optimization approach},
   volume = {239},
   url = {http://dx.doi.org/10.1016/j.physd.2010.04.001},
   year = {2010},
}

@article{Tu2010,
   author = {Z. C. Tu},
   doi = {10.1063/1.3335894},
   issn = {0021-9606},
   issue = {8},
   journal = {The Journal of Chemical Physics},
   month = {2},
   pmid = {20192294},
   title = {Compatibility between shape equation and boundary conditions of lipid membranes with free edges},
   volume = {132},
   url = {https://pubs.aip.org/jcp/article/132/8/084111/188738/Compatibility-between-shape-equation-and-boundary},
   year = {2010},
}

@article{Arroyo2009,
   author = {Marino Arroyo and Antonio Desimone},
   doi = {10.1103/PhysRevE.79.031915},
   issn = {15393755},
   issue = {3},
   journal = {Physical Review E - Statistical, Nonlinear, and Soft Matter Physics},
   pages = {1-17},
   title = {Relaxation dynamics of fluid membranes},
   volume = {79},
   year = {2009},
}

@article{Srividya2008,
   author = {Narayanan Srividya and Subra Muralidharan},
   doi = {10.1021/jp7119203},
   issn = {1520-6106},
   issue = {24},
   journal = {The Journal of Physical Chemistry B},
   month = {6},
   pages = {7147-7152},
   pmid = {18503265},
   title = {Determination of the Line Tension of Giant Vesicles from Pore-Closing Dynamics},
   volume = {112},
   url = {https://pubs.acs.org/doi/10.1021/jp7119203},
   year = {2008},
}

@article{Alexander2006,
   author = {James C. Alexander and Andrew J. Bernoff and Elizabeth K. Mann and J. Adin Mann and Lu Zou},
   doi = {10.1063/1.2212887},
   issn = {1070-6631},
   issue = {6},
   journal = {Physics of Fluids},
   month = {6},
   title = {Hole dynamics in polymer Langmuir films},
   volume = {18},
   url = {https://pubs.aip.org/pof/article/18/6/062103/919710/Hole-dynamics-in-polymer-Langmuir-films},
   year = {2006},
}

@article{Umeda2005,
   author = {Tamiki Umeda and Yukio Suezaki and Kingo Takiguchi and Hirokazu Hotani},
   doi = {10.1103/PhysRevE.71.011913},
   issn = {1539-3755},
   issue = {1},
   journal = {Physical Review E},
   month = {1},
   pages = {011913},
   pmid = {15697636},
   title = {Theoretical analysis of opening-up vesicles with single and two holes},
   volume = {71},
   url = {https://link.aps.org/doi/10.1103/PhysRevE.71.011913},
   year = {2005},
}

@article{Yin2005,
   author = {Yajun Yin and Jie Yin and Dong Ni},
   doi = {10.1007/s00285-005-0330-x},
   issn = {03036812},
   issue = {4},
   journal = {Journal of Mathematical Biology},
   pages = {403-413},
   pmid = {15940540},
   title = {General mathematical frame for open or closed biomembranes (Part I): Equilibrium theory and geometrically constraint equation},
   volume = {51},
   year = {2005},
}

@article{Yin2005b,
   author = {Yajun Yin and Yanqiu Chen and Dong Ni and Huiji Shi and Qinshan Fan},
   doi = {10.1016/j.jbiomech.2004.06.024},
   issn = {00219290},
   issue = {7},
   journal = {Journal of Biomechanics},
   pages = {1433-1440},
   pmid = {15922754},
   title = {Shape equations and curvature bifurcations induced by inhomogeneous rigidities in cell membranes},
   volume = {38},
   year = {2005},
}

@article{Jiang2004,
   author = {Shidong Jiang and Vladimir Rokhlin},
   doi = {10.1016/j.jcp.2003.10.001},
   issn = {00219991},
   issue = {1},
   journal = {Journal of Computational Physics},
   pages = {1-16},
   title = {Second kind integral equations for the classical potential theory on open surfaces II},
   volume = {195},
   year = {2004},
}

@article{Tu2003,
   author = {Z. C. Tu and Z. C. Ou-Yang},
   doi = {10.1103/PhysRevE.68.061915},
   issn = {1063651X},
   issue = {6},
   journal = {Physical Review E - Statistical Physics, Plasmas, Fluids, and Related Interdisciplinary Topics},
   pages = {1-7},
   pmid = {14754242},
   title = {Lipid membranes with free edges},
   volume = {68},
   year = {2003},
}

@article{Jiang2003,
   author = {Shidong Jiang and Vladimir Rokhlin},
   doi = {10.1016/S0021-9991(03)00304-8},
   issn = {00219991},
   issue = {1},
   journal = {Journal of Computational Physics},
   pages = {40-74},
   title = {Second kind integral equations for the classical potential theory on open surfaces I: Analytical apparatus},
   volume = {191},
   year = {2003},
}

@article{Misbah2006,
   author = {Chaouqi Misbah},
   doi = {10.1103/PhysRevLett.96.028104},
   issn = {0031-9007},
   issue = {2},
   journal = {Physical Review Letters},
   month = {1},
   pages = {028104},
   pmid = {16486649},
   title = {Vacillating Breathing and Tumbling of Vesicles under Shear Flow},
   volume = {96},
   url = {https://link.aps.org/doi/10.1103/PhysRevLett.96.028104},
   year = {2006},
}

@article{Lai2001,
   author = {Ming Chih Lai and Zhilin Li},
   doi = {10.1016/S0893-9659(00)00127-0},
   issn = {08939659},
   issue = {2},
   journal = {Applied Mathematics Letters},
   pages = {149-154},
   title = {A remark on jump conditions for the three-dimensional Navier-Stokes equations involving an immersed moving membrane},
   volume = {14},
   year = {2001},
}

@article{Cortez2005,
   author = {Ricardo Cortez and Charles S. Peskin and John M. Stockie and Douglas Varela},
   doi = {10.1137/S003613990342534X},
   issn = {00361399},
   issue = {2},
   journal = {SIAM Journal on Applied Mathematics},
   pages = {494-520},
   title = {Parametric resonance in immersed elastic boundaries},
   volume = {65},
   year = {2005},
}

@article{Lai2022,
   author = {Ming Chih Lai and Yunchang Seol},
   doi = {10.1016/j.jcp.2021.110785},
   issn = {10902716},
   journal = {Journal of Computational Physics},
   month = {1},
   pages = {110785},
   publisher = {Elsevier Inc.},
   title = {A stable and accurate immersed boundary method for simulating vesicle dynamics via spherical harmonics},
   volume = {449},
   url = {https://doi.org/10.1016/j.jcp.2021.110785 https://linkinghub.elsevier.com/retrieve/pii/S002199912100680X},
   year = {2022},
}

@article{HSU2019747,
   author = {Shih Hsuan Hsu and Jay Chu and Ming Chih Lai and Richard Tsai},
   doi = {10.1016/j.jcp.2019.06.046},
   issn = {10902716},
   journal = {Journal of Computational Physics},
   pages = {747-764},
   publisher = {Elsevier Inc.},
   title = {A coupled grid based particle and implicit boundary integral method for two-phase flows with insoluble surfactant},
   volume = {395},
   url = {https://doi.org/10.1016/j.jcp.2019.06.046 https://www.sciencedirect.com/science/article/pii/S0021999119304541},
   year = {2019},
}

@misc{Helsing2024,
      title={The Helmholtz Dirichlet and Neumann problems on piecewise smooth open curves}, 
      author={Johan Helsing and Shidong Jiang},
      year={2024},
      eprint={2411.05761},
      archivePrefix={arXiv},
      primaryClass={math.NA},
      url={https://arxiv.org/abs/2411.05761}, 
}

@article{Hayashi1977,
   author = {Yoshio Hayashi},
   doi = {10.3792/pjaa.53.159},
   isbn = {2013206534},
   issn = {0386-2194},
   issue = {5},
   journal = {Proceedings of the Japan Academy, Series A, Mathematical Sciences},
   month = {10},
   pages = {1-23},
   title = {Three-dimensional Dirichlet problem for the Helmholtz equation for an open boundary},
   volume = {53},
   url = {https://projecteuclid.org/journals/proceedings-of-the-japan-academy-series-a-mathematical-sciences/volume-53/issue-5/Three-dimensional-Dirichlet-problem-for-the-Helmholtz-equation-for-an/10.3792/pjaa.53.159.full},
   year = {1977},
}

@article{Kirvalidze1997,
   author = {V. Kirvalidze},
   doi = {10.1002/(SICI)1099-1476(199710)20:15<1257::AID-MMA896>3.0.CO;2-5},
   issn = {0170-4214},
   issue = {15},
   journal = {Mathematical Methods in the Applied Sciences},
   month = {10},
   pages = {1257-1269},
   title = {The Dirichlet problem for Stokes equation in a domain exterior to an open surface},
   volume = {20},
   url = {https://onlinelibrary.wiley.com/doi/10.1002/(SICI)1099-1476(199710)20:15%3C1257::AID-MMA896%3E3.0.CO;2-5},
   year = {1997},
}

@article{Stephan1987,
   author = {Ernst P. Stephan},
   doi = {10.1007/BF01199079},
   issn = {0378-620X},
   issue = {2},
   journal = {Integral Equations and Operator Theory},
   month = {3},
   pages = {236-257},
   title = {Boundary integral equations for screen problems in IR3},
   volume = {10},
   url = {http://link.springer.com/10.1007/BF01199079},
   year = {1987},
}

@article{Bottacchiari2024,
author = {Bottacchiari, Matteo and Gallo, Mirko and Bussoletti, Marco and Casciola, Carlo M.},
doi = {10.1093/pnasnexus/pgae300},
editor = {Amon, Cristina},
issn = {2752-6542},
journal = {PNAS Nexus},
month = {aug},
number = {8},
pages = {1--10},
publisher = {Oxford University Press},
title = {{The diffuse interface description of fluid lipid membranes captures key features of the hemifusion pathway and lateral stress profile}},
url = {https://doi.org/10.1093/pnasnexus/pgae300 https://academic.oup.com/pnasnexus/article/doi/10.1093/pnasnexus/pgae300/7720659},
volume = {3},
year = {2024}
}

@article{Wang2007,
author = {Wang, Xiaoqiang and Du, Qiang},
doi = {10.1007/s00285-007-0118-2},
issn = {0303-6812},
journal = {Journal of Mathematical Biology},
month = {nov},
number = {3},
pages = {347--371},
pmid = {17701177},
title = {{Modelling and simulations of multi-component lipid membranes and open membranes via diffuse interface approaches}},
url = {http://link.springer.com/10.1007/s00285-007-0118-2},
volume = {56},
year = {2007}
}

@article{Cohen2012,
author = {Cohen, Fredric S. and Eisenberg, Robert and Ryham, Rolf J.},
doi = {10.4310/CMS.2012.v10.n4.a12},
issn = {15396746},
journal = {Communications in Mathematical Sciences},
number = {4},
pages = {1273--1285},
title = {{A dynamic model of open vesicles in fluids}},
url = {https://link.intlpress.com/JDetail/1806264961234063361},
volume = {10},
year = {2012}
}

@article{Cohen2014,
author = {Cohen, Fredric S. and Ryham, Rolf J.},
doi = {10.1063/1.4864192},
issn = {1070-6631},
journal = {Physics of Fluids},
month = {feb},
number = {2},
title = {{The aqueous viscous drag of a contracting open surface}},
url = {https://pubs.aip.org/pof/article/26/2/023101/259225/The-aqueous-viscous-drag-of-a-contracting-open},
volume = {26},
year = {2014}
}

@book{mclean2000strongly,
  author    = {McLean, William},
  title     = {Strongly Elliptic Systems and Boundary Integral Equations},
  publisher = {Cambridge University Press},
  address   = {Cambridge},
  year      = {2000},
  series    = {Cambridge Monographs on Applied and Computational Mathematics},
  volume    = {15}
}

@article{VonPetersdorff1990,
author = {{Von Petersdorff}, T. and Stephan, E. P.},
doi = {10.1002/mma.1670120306},
issn = {10991476},
journal = {Mathematical Methods in the Applied Sciences},
number = {3},
pages = {229--249},
title = {{Regularity of mixed boundary value problems in $\mathbb R^3$ and boundary element methods on graded meshes}},
volume = {12},
year = {1990}
}

@article{kuchta2020mixed,
  author    = {Miroslav Kuchta and Federica Laurino and Kent-Andre Mardal and Paolo Zunino},
  title     = {Analysis and approximation of mixed-dimensional PDEs on 3D–1D domains coupled with Lagrange multipliers},
  journal   = {SIAM Journal on Numerical Analysis},
  volume    = {58},
  number    = {6},
  pages     = {3673--3700},
  year      = {2020},
  doi       = {10.1137/20M1329664}
}

@article{grappein2024xfem,
  author    = {Denise Grappein and Stefano Scialò and Fabio Vicini},
  title     = {Extended finite elements for 3D–1D coupled problems via a PDE-constrained optimization approach},
  journal   = {Computer Methods in Applied Mechanics and Engineering},
  volume    = {415},
  pages     = {116237},
  year      = {2024},
  doi       = {10.1016/j.cma.2024.116237}
}

@article{quarteroni2007stability,
  author    = {Alfio Quarteroni and Alessandro Veneziani},
  title     = {On the stability of the coupling of 3D and 1D fluid-structure interaction models},
  journal   = {ESAIM: Mathematical Modelling and Numerical Analysis},
  volume    = {41},
  number    = {4},
  pages     = {697--721},
  year      = {2007},
  doi       = {10.1051/m2an:2007030}
}

@article{Boedec2017,
author = {Boedec, Gwenn and Leonetti, Marc and Jaeger, Marc},
issn = {0021-9991},
journal = {Journal of Computational Physics},
pages = {117--138},
publisher = {Elsevier Inc.},
title = {{Isogeometric FEM-BEM simulations of drop , capsule and vesicle dynamics in Stokes flow}},
volume = {342},
year = {2017}
}

@article{Adkins2025,
author = {Adkins, Raymond and Robaszewski, Joanna and Shin, Seungwoo and Brauns, Fridtjof and Jia, Leroy and Khanra, Ayantika and Sharma, Prerna and Pelcovits, Robert A and Powers, Thomas R and Dogic, Zvonimir},
doi = {10.1073/pnas.2427024122},
journal = {Proceedings of the National Academy of Sciences},
number = {36},
pages = {e2427024122},
title = {{Topology and kinetic pathways of colloidosome assembly and disassembly}},
url = {https://www.pnas.org/doi/abs/10.1073/pnas.2427024122},
volume = {122},
year = {2025}
}
